\documentclass[11pt]{amsart}
\pdfoutput=1
\usepackage{amsmath, amsthm, amssymb}
\usepackage{graphicx}
\usepackage{mathbbol}
\usepackage{txfonts}
\usepackage[usenames]{color}
\usepackage{srcltx} 

\definecolor{violet}{RGB}{100,0,100}

\newtheorem{thm}{Theorem}[section]
\newcommand{\bt}{\begin{thm}}
\newcommand{\et}{\end{thm}}

\newtheorem{conj}[thm]{Conjecture}

\newtheorem{question}[thm]{Question}

\newtheorem{quest}[thm]{Question}

\newtheorem{cor}[thm]{Corollary}   
\newcommand{\bc}{\begin{cor}}
\newcommand{\ec}{\end{cor}}

\newtheorem{lem}[thm]{Lemma}   
\newcommand{\bl}{\begin{lem}}
\newcommand{\el}{\end{lem}}

\newtheorem{prop}[thm]{Proposition}
\newcommand{\bp}{\begin{prop}}
\newcommand{\ep}{\end{prop}}

\newtheorem{defn}[thm]{Definition}
\newcommand{\bd}{\begin{defn}}    
\newcommand{\ed}{\end{defn}}

\newtheorem{rmrk}[thm]{Remark}   

\newcommand{\br}{\begin{rmrk}}
\newcommand{\er}{\end{rmrk}}

\newcommand{\GHto}{\xrightarrow{\textrm{\,GH\,\,}} }
\newcommand{\FFto}{\xrightarrow{\textrm{\, F\,\,}} }
\newcommand{\Fto}{\xrightarrow{\mathcal{\,F\,\,}} }
\newcommand{\mto}{\xrightarrow{ {\,\,m\,\,}} }
\newcommand{\VFto}{\xrightarrow{\mathcal{\,VF\,\,}} }
\newcommand{\VADBto}{\xrightarrow{\textrm{\,VADB\,\,}} }
\newcommand{\LIPto}{\xrightarrow{\textrm{\,Lip\,\,}} }

\newcommand{\be}{\begin{equation}}

\newcommand{\ee}{\end{equation}}

\newcommand{\Width}{\operatorname{Width}}
\newcommand{\Diam}{\operatorname{Diam}}

\newcommand{\Vol}{\operatorname{Vol}}
\newcommand{\FillVol}{\operatorname{FillVol}}
\newcommand{\Scal}{{\rm Scal}} 
\newcommand{\MinA}{\operatorname{MinA}}
\newcommand{\Area}{\operatorname{Area}}

\newcommand{\Ricci}{\rm{Ricci}}

\newcommand{\disjointunion}{\sqcup}

\newcommand{\Lip}{\operatorname{Lip}}

\newcommand{\mass}{{\mathbf M}}

\newcommand{\Length}{\operatorname{Length}}

\newcommand{\vol}{\operatorname{Vol}}

\begin{document}

\title{Conjectures on Convergence and Scalar Curvature}

\author[Sormani et. al.]{Christina Sormani and Participants at the IAS Emerging Topics Workshop on Scalar Curvature and Convergence}
 \thanks{The Institute for Advanced Study and NSF DMS 1006059.}
\address{CUNY Graduate Center and Lehman College}
\email{sormanic@gmail.com}





\begin{abstract} 
Here we survey the compactness and geometric stability conjectures formulated by the participants at the 2018 IAS Emerging Topics Workshop on {\em Scalar Curvature and Convergence}.   We have tried to survey all the progress towards these conjectures as well as related examples, although it is impossible to cover everything.  We focus primarily on sequences of compact Riemannian manifolds with nonnegative scalar curvature and their limit spaces.  Christina Sormani is grateful to have had the opportunity to write up our ideas and has done her best to credit everyone involved within the paper even though she is the only author listed above.  In truth we are a team of over thirty people working together and apart on these deep questions and we welcome everyone who is interested in these conjectures to join us.
\end{abstract}

\maketitle

\newpage

\section{Introduction}

\begin{figure}[h]
\begin{center}
\includegraphics[width=0.8\textwidth]{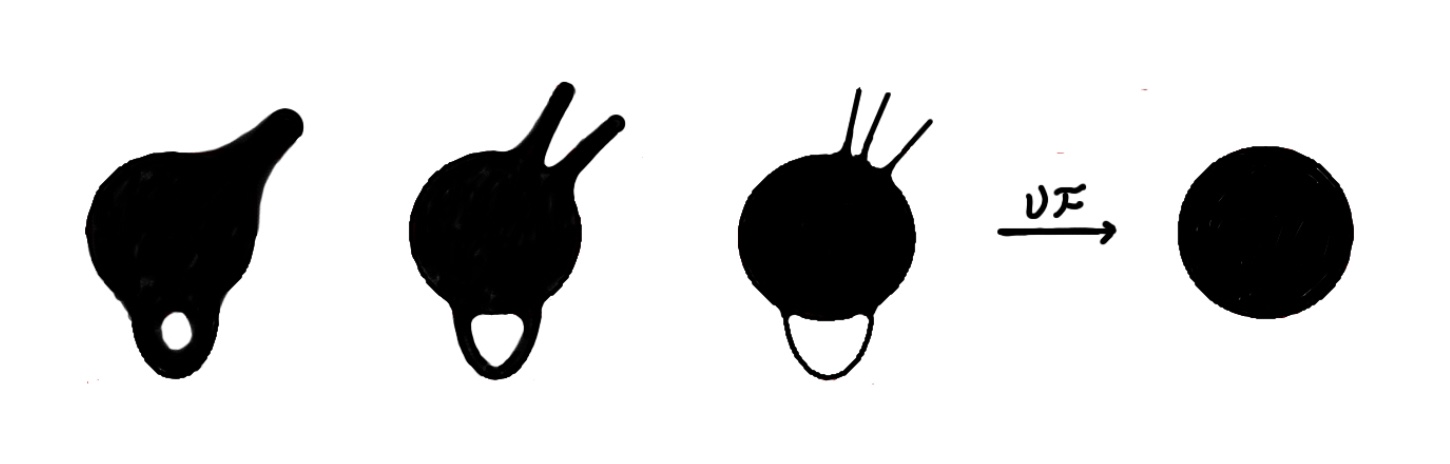} 
\caption{A sequence of $M_j^3$ with positive scalar curvature that converges to a round sphere, ${\mathbb S}^3$, in the $\mathcal{VF}$ sense but have no GH limit.}\label{VF-Example}
\end{center}
\end{figure}

One of the greatest challenges in Geometric Analysis today is to develop a deeper understanding of the 
geometry of manifolds with lower bounds on their scalar curvature and their limit spaces. Even in three dimensions, manifolds with positive scalar curvature can have arbitrarily thin and deep wells and 
arbitrarily long or short tunnels (see Figure~\ref{VF-Example}).   The existence of a single increasingly thin well prevents the existence of a smooth or Lipschitz limit.   If there are increasingly many wells, the sequence has no 
Gromov-Hausdorff (GH) limit either.  For this reason we must consider weaker notions of convergence and at the same time work to develop a notion of convergence which is strong enough to preserve the key properties of scalar curvature. 
In Section~\ref{Sect-scalar} we review the key geometric properties of scalar curvature including new low regularity notions of nonnegative scalar curvature. 

In \cite{Gromov-Dirac}, Gromov suggested that we apply the notion of intrinsic flat ($\mathcal{F}$) convergence 
which was defined by Sormani-Wenger in \cite{SW-JDG} for sequences of oriented Riemannian manifolds with boundary using Ambrosio-Kirchheim theory \cite{AK}.    Intrinsic flat convergence has the advantage of having an existing compactness theorem by Wenger which requires only global volume and diameter bounds while at the same time producing rectifiable limit spaces \cite{Wenger-compactness} \cite{SW-JDG}.   Although points may disappear in the limit, balls about points which do not disappear will converge to balls about the limit points, the boundaries of those balls will converge, and their filling volumes will converge.     In fact volume preserving intrinsic flat ($\mathcal{VF}$) convergence, (which is $\mathcal{F}$ convergence combined with a requirement that the total volumes converge) even implies convergence in measure.   There have already been significant advances exploring how $\mathcal{VF}$ convergence can be applied to understand the limits of sequences of Riemannian manifolds with uniform lower bounds on their scalar curvature.   In Section~\ref{SWIF}, we will survey the definitions and properties of $\mathcal{F}$, $\mathcal{VF}$, and other notions of weak convergence that might be applied to explore the geometry of Riemannian manifolds with scalar curvature bounds.

Here we present seven key conjectures that were formulated or reformulated by various mathematicians during the {\em 2018 Emerging Topics at IAS on Scalar Curvature and Convergence}.  
There are three compactness conjectures:
\begin{itemize}
\item{MinA Scalar Compactness [Gromov-Sormani]}
\item{NoMin Boundary Scalar Compactness [Gromov]}
\item{Scalar Mass Compactness [Sormani]}
\end{itemize}
each of which has different suggested hypotheses that would provide enough control on a sequence of manifolds, $M_j$, with nonnegative scalar curvature so that a subsequence $\mathcal{VF}$ converges to a limit space which has some generalized notion of nonnegative scalar curvature.   There are four geometric stability conjectures:
\begin{itemize}
\item{Geometric Stability of Scalar Torus Rigidity [Gromov-Sormani]}
\item{Geometric Stability of Scalar Prism Rigidity [Gromov-Li]}
\item{Geometric Stability of Scalar Sphere Rigidity [Marques-Neves]}
\item{Geometric Stability of Zero Mass Rigidity [Lee-Sormani]} 
\end{itemize}
each of which states that a manifold which almost satisfies the hypotheses of a known rigidity theorem must be 
$\mathcal{VF}$ close to the rigid manifold that exactly satisfies those hypotheses. If we prove one of the compactness conjectures and prove the rigidity theorems hold on the class of limit spaces, then these geometric stability conjectures would follow.  They might also be proven directly without a compactness theorem either in full generality or in special settings.   All these conjectures concern Riemannian manifolds which are three dimensional.

Each of these seven conjectures will be discussed in its own section where we will provide a precise statement of the conjecture, review the history of the conjecture, survey progress towards the conjecture,  discuss related examples, and finally describe the consequences of each conjecture.  Throughout we will also list open problems at various levels of difficulty.  Some of the open problems are suggested special cases of the conjectures while others involve proving the consequences of the conjectures directly.   

In the final section of the paper we return to the properties of $\mathcal{VF}$ and VADB convergence with brief survey of known results and a list of open questions.  As it impossible to discuss all conjectures and ideas that arose during discussions at the {\em Emerging Topics at IAS on Scalar Curvature and Convergence} (IAS) and subsequent workshops including the {\em 2020 Virtual Workshop on Ricci and Scalar Curvature} (VWRS) we recommend the
reader consult Gromov's surveys in \cite{Gromov-Natural} and \cite{Gromov-Four-Lectures} and other chapters within this volume.

\vspace{.3cm}
\noindent{\bf Acknowledgements:} The writer would like to thank all the participants in the workshop at IAS and subsequent events including Misha Gromov, Brian Allen, Lucas Ambrozio, Alessandro Carlotto, Otis Chodosh, Fernando Coda Marques, Michael Eichmair, Bernhard Hanke, Lan-Hsuan Huang, Jeff Jauregui, Nicos Kapouleas, Demetre Kazaras, Anusha Krishnan, Sajjad Lakzian, Dan Lee, Chao Li, Yevgeny Liokumovich, Siyuan Lu, Elena Maeder-Baumdicker, Andrea Malchiodi, Yashar Memarian, Pengzi Miao, Alex Nabutovsky, Robin Neumayer, Andre Neves, Raquel Perales, Jacobus Portegies,  Regina Rotman, Rick Schoen, Shengwen Wang, Guofang Wei, Ruobing Zhang, and Xin Zhou for deep conversations at IAS.    We would also like to thank Yuguang Shi, Robert Young, Armando Cabrera Pacheco, Jorge Basilio, Lisandra Hernandez Vasquez, 
Paula Burkhardt-Guim, Daniel Stern, Edward Bryden, Marcus Khuri,
James Isenberg, Martin Lesourd, Chritian Ketterer, Armando Cabrera Pacheco, Aaron Naber, Man-Chun Lee, Robin Neumeyer, Nicola Gigli, Shoehei Honda, Richard Bamler, Jintian Zhu, Guofang Wei, Paolo Piazza, and Blaine Lawson for conversations at Stony Brook, NYU, and VWRS.   We 
try to credit everyone within for their specific suggestions where possible.   I would especially like to give thanks to Raquel Perales, Jeff Jauregui, Chao Li, Pengzi Miao, Dan Lee, Brian Allen,  Lan-Hsuan Huang, and Michael Eichmair for extra assistance in the final stages of preparation of this work.

\section{The Geometry of Scalar Curvature} \label{Sect-scalar}

Before we state our conjectures or define any notions of convergence, we review some of the key geometric properties of three dimensional Riemannian manifolds with scalar curvature bounds.  We have chosen to focus on the following four rigidity theorems: 
\begin{itemize}
\item{The Scalar Torus Rigidity Theorem [Schoen-Yau] [Gromov-Lawson]}
\item{The Scalar Prism Rigidity Theorem [Gromov] [Li]}
\item{The Scalar Sphere Rigidity Theorem [Marques-Neves]}
\item{Zero Mass Rigidity [Schoen-Yau][Witten][Shi-Tam]} 
\end{itemize}
Note that we quite deliberately avoid using differentiability in the definitions and in the statements of these theorems.  This allows us to see how these notions might be defined and theorems might be stated on limit spaces of lower regularity.

Let us begin by recalling that a compact Riemannian manifold, $(M, g)$, with a $C^0$ metric tensor $g$ may be viewed as a compact metric space, $(M,d)$, by taking
\be \label{d-from-g}
d(p,q) = \inf\{ L_g(C)\, |\, C(0)=p,\, C(1)=q\} \textrm{ where } L(C)=\int_0^1 g(C'(t), C'(t))^{1/2} \, dt.
\ee
This infimum is achieved by a geodesic, so we say $(M,d)$ is a {\em compact geodesic metric space}. 
The diameter of the manifold is 
\be
\Diam(M)= \sup\{ d(x,y)\,: \, x,y\in M\}.
\ee
Balls are then defined $B(p,r)=\{x\,|\, d(x,p)<r\}$ and the exponential map provides a biLipschitz chart from
a ball in Euclidean space to a ball in $(M,d)$.   The manifold is oriented if a preferred orientation has been chosen for all charts so that that Jacobian of the biLipschitz transition maps has positive determinant almost everywhere.    

Volumes of regions, $U\subset M$, can be defined either by integrating $\sqrt{\det g}$ on disjoint charts or equivalently using the Hausdorff measure
\be
\Vol(U)= \mathcal{H}^n(U) 
= \lim_{\delta\to 0}  \inf\left\{c_n \sum_{j=1}^\infty (\Diam(U_j))^n\, : \, U \subset \bigcup_{j=1}^\infty U_j, \,\, \Diam(U_j) \le \delta\right\}.
\ee

To compute the scalar curvature at $p$ we need only take the limit
\be \label{Scal-vol-lim}
\Scal(p)=\lim_{r\to 0} 6(n+2) \left(\frac{V^n(r) - \Vol(B(p,r))}{r^2 V^n(r)} \right)
\ee
where $V^n(r)$ is the volume of a Euclidean ball of radius $r$ in dimension $n$.
Note that a {\em three dimensional round sphere}, $({\mathbb{S}}^3, g_{{\mathbb{S}}^3})$, with diameter $\pi$ has $\Scal=6$ everywhere.   An
example with $\Scal=0$ everywhere is a
 {\em flat three torus}, $(M^3, g_0)$ which is a manifold homeomorphic to a torus 
 \be
 {\mathbb{T}}^3={\mathbb{S}}^1\times {\mathbb{S}}^1\times{\mathbb{S}}^1
 \ee
  such that balls about every point are isometric to Euclidean balls.   

\subsection{Stable Minimal Surfaces}

\begin{figure}[h]
\begin{center}
\includegraphics[width=0.9\textwidth]{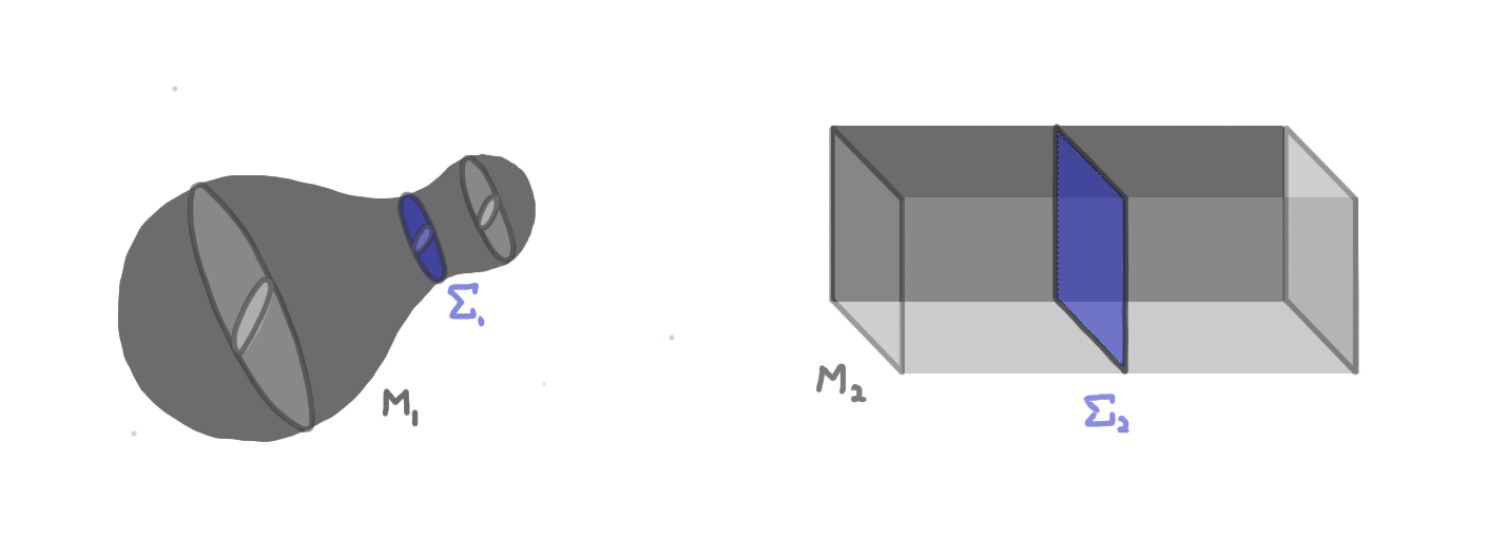} 
\caption{Two manifolds $M_i$ with $\Scal \ge 0$: the first contains a stable minimal sphere, 
$\textcolor{blue}{\Sigma_1}\subset M_1$,
and the second has a stable minimal torus $\textcolor{blue}{\Sigma_2}\subset M_2$.  Thus $M_2$ is isometric to a flat three torus.}\label{fig-torus-rigidity}
\end{center}
\end{figure}

The area of a surface may also be defined either by integration or equivalently by using the Hausdorff measure.
A surface $\Sigma^2\subset M^3$ is a closed minimal surface if it has no boundary and locally minimizes area:
\be
\Area(\Sigma^2\cap U) \le \Area(\Sigma') \qquad \forall \Sigma' \,\,s.t.\,\, \partial \Sigma'=\partial(\Sigma^2\cap U).
\ee
It is a stable closed minimal surface if it minimizes area of any continuous deformation, $\Sigma_t$:
\be
\Area(\Sigma^2) \le \Area(\Sigma_t) \qquad \forall \Sigma_t \,\,s.t.\,\, \Sigma_0=\Sigma.
\ee
In \cite{Schoen-Yau-minimal}, Richard Schoen and Shing-Tung Yau proved that if $g$ is $C^2$ then by taking the second variation of the area and applying the Gauss-Bonnet Theorem, the only stable closed minimal surfaces in a manifold with nonnegative scalar curvature are tori and spheres.   If a compact $M^3$ with nonnegative scalar curvature contains a minimal  torus, then $M^3$ is a flat three dimensional torus.   See Figure~\ref{fig-torus-rigidity}. Thus in particular they proved the following:

\begin{thm} {\bf Scalar Torus Rigidity  [Schoen-Yau]}
If $M^3$ is homeomorphic to a torus and has nonnegative scalar curvature then is isometric to a flat torus.
\end{thm}.

This theorem was generalized to higher dimensions by Misha Gromov and Blaine Lawson \cite{Gromov-Lawson-Spin1} with a new proof that also works for dimension three.   Recently Stern has presented a new proof using harmonic maps from the torus to the circle ${\mathbb{S}}^1$ \cite{Stern-Scalar-Harmonic}.   

In addition to the Scalar Torus Rigidity Theorem, there are a number of other rigidity theorems which similarly involve the existence of a stable minimal or CMC surface in a Riemannian manifold with lower bounds on scalar curvature.   
In the asymptotically flat setting there is the Rigidity of the Positive Mass Theorem of Schoen-Yau  \cite{Schoen-Yau-PMT}
which we will discuss later along with the hyperbolic version of this theorem and other noncompact theorems.
In the compact setting there is the
Minimal Torus Scalar Rigidity Theorem of Cai-Galloway \cite{CG-CAG2000}, the
Hemispherical Scalar Rigidity Theorem of Eichmair \cite{Eichmair-Hemi09}, the Cover Splitting Scalar Rigidity Theorem of Bray, Brendle, and Neves \cite{BBN10},  and the ${\mathbb{RP}}^3$
Scalar Rigidity Theorem by Bray, Brendle, Eichmair and Neves \cite{BBEN10}.  All of these involve minimal surfaces in three dimensional manifolds, but there is also one involving minimal surfaces in four dimensional manifolds by Zhu
in \cite{Zhu-2inNdim}.   Ambrozio has a rigidity theorem involving stable minimal surfaces with free boundary \cite{Ambrozio-free}.   

We should also keep in mind the lack of rigidity 
discovered in work of Corvino \cite{Corvino-deform}, Lohkamp \cite{Lohkamp-hammocks}, Hang-Wang \cite{Hang-Wang-non}, and Brendle-Marques-Neves \cite{BMN-False-Oo}.  We should also keep in mind the
Schoen-Yau and Gromov-Lawson tunnel constructions in \cite{Schoen-Yau-structure} and \cite{Gromov-Lawson-classification}, and the more recent sewing constructions of Basilio-Dodziuk-Sormani in \cite{BDS-sewing} using those tunnels.  In dimension 4 and up there is an intriguing set of examples by Lee-Naber-Neumayer in \cite{Lee-Naber-Neumayer-dp} using a new method altogether.

\subsection{Prism Rigidity of Gromov and Li}

Gromov proposed the following theorem in \cite{Gromov-Dirac} and it was proven using Schoen-Yau style methods
by Chao Li 
in \cite{Li-polyhedral}.  For simplicity we state the theorem for cubes but Li's proof holds for more general prisms. 

\begin{thm} {\bf Prism Rigidity  [Gromov][Li]}  
If $M^3$ has $\Scal \ge 0$, then any cube $P \subset M^3$ satisfies prism rigidity.  That is:
if the faces of $P^3$ are mean convex and dihedral angles, $\theta_p \le \pi/2$, at any $p$ lying on 
an edge of $P^3$,
then $P$ is isometric to a rectangular prism in Euclidean space.
\end{thm}

\begin{figure}[h]
\begin{center}
\includegraphics[width=0.8\textwidth]{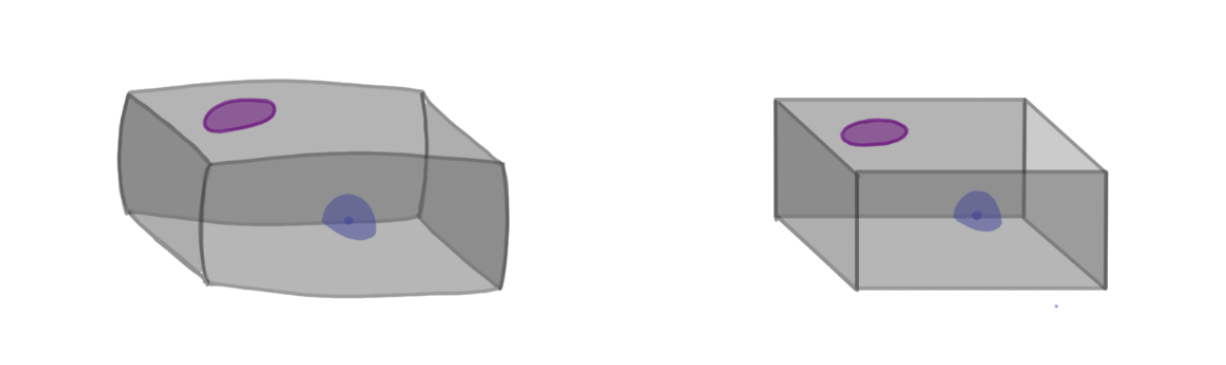} 
\caption{On the left is a prism $P\subset M^3$ with  mean convex boundaries
has \textcolor{violet}{$\Area(\Sigma_t)\le \Area(\Sigma_0)$} for local deformations
that has dihedral angles \textcolor{blue}{$\theta_{p} \le  \pi/2$}.   If 
$\Scal\ge 0$ on $P$ then $P$ is isometric to a rectangular prism in Euclidean space as depicted on the right.}\label{fig-prism}
\end{center}
\end{figure}

Mean convexity of the boundaries can be defined by stating that local continuous inward deformations of the boundary
are area nonincreasing:
\be
\Area(\Sigma_t) \le \Area(\Sigma_0)
\ee
where $\Sigma_0\subset \partial M$ lies in a face of $\partial M$ avoiding an edge, and $\Sigma_t \subset M$ is a variation of $\Sigma_0$ with $\partial \Sigma_t=\partial \Sigma_0$.  So this includes boundaries whose faces are minimal surfaces.  

The dihedral angle at $p\in \partial P$ can be defined using volumes:
\be
\theta_p(P) = \lim_{r\to 0} 2 \pi \vol(B(p,r)\cap P) / V^3(r).
\ee
 See Figure~\ref{fig-prism}.  Note that the hypotheses of the theorem only depend on area and volume.
 In the conclusion distances are controlled because the lengths, $x$, $y$, and $z$, of the sides of a 
 rectangular prism in Euclidean space are determined by the total volume divided by area of the perpendicular face: $x=(xyz)/(yz)$.  

Gromov's proof of prism rigidity of cubes begins by reflecting the cube across the mean convex faces to create a torus with a $C^0$ metric tensor.  He then smooths the metric on the cube to a metric with nonnegative scalar curvature
to reduce the problem to the torus rigidity theorem \cite{Gromov-Dirac}.   He also proves that if all the small cubes in a smooth manifold $M^3$  satisfy the prism rigidity property, then
$M^3$ has nonnegative Scalar curvature.   He has proposed this as a natural low regularity definition of nonnegative scalar curvature.

\begin{quest} If a $C^0$ manifold satisfies prism rigidity for small prisms, is the limit of the ratio of volumes in (\ref{Scal-vol-lim}) nonnegative?
\end{quest}

Chao Li proves the prism rigidity theorem for more general prisms using a capillary flow \cite{Li-polyhedral}.   His prisms do not need to be cubes.  The faces are still mean convex and the dihedral angles are chosen in comparison to a corresponding Euclidean prism.   It is important to note there are no assumptions on the lengths of edges or the areas of faces in these prism rigidity theorems.   This version of the proof does not apply the Torus Rigidity Theorem and is instead a direct proof which might be said to be in the style of a Schoen-Yau minimal surface proof as the levels of the capillary flow are minimal surfaces, but the boundaries are not free and they have corners.  He applies  Jean Taylor's proof of
the existence of capillary surfaces as integral currents \cite{Taylor-capillary} and then improves the regularity of the surfaces so that he can apply the technique Bray-Brendle-Neves used in \cite{BBN10} to prove the Scalar Cover Splitting Theorem.   He has also proven a hyperbolic version of this theorem assuming $\Scal \ge -1$ in \cite{Li-dihedral-hyp}.

\subsection{Widths and Unstable Minimal Surfaces}

\begin{figure}[h]
\begin{center}
\includegraphics[width=0.6\textwidth]{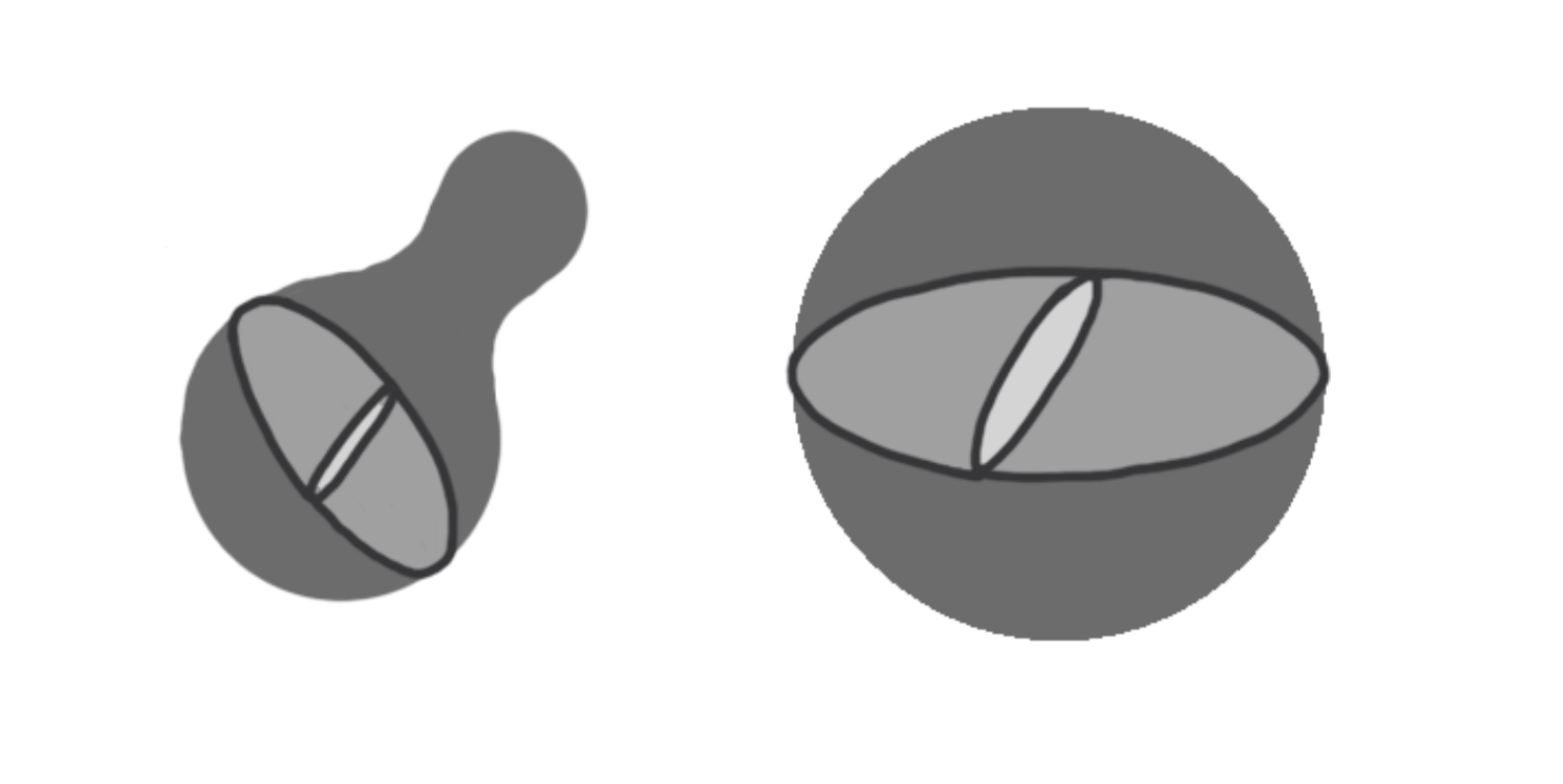} 
\caption{The Riemannian manifold on the left has $\Scal \ge 6$ everywhere and an unstable minimal surface achieving the width of area $< 4\pi$ while the manifold on the right has $\Scal \ge 6$ everywhere and an unstable minimal surface achieving the width of area $=4\pi$ so it is isometric to ${\mathbb S}^3$.}\label{fig-sphere}
\end{center}
\end{figure}

Marques and Neves have proven a rigidity theorem involving an unstable minimal surface in \cite{MN-Duke}.   Let us begin by reviewing the definition of Width:
\be
\Width(M^3) = \inf_{\Sigma_t\in \Lambda} \sup_{t\in [0,1]} \Area(\Sigma_t)
\ee
where $\Lambda$ is the collection of all sweepouts of $M^3$.   More precisely, when $M^3$ is a sphere, viewed
as a subset of ${\mathbb{E}}^4$, a sweepout consists of a foliation by surfaces:
\be
\Sigma_t= F_t\left(x_4^{-1}(t)\right) \textrm{ where } F_t: M^3 \to M^3
\ee
is a smooth family of diffeomorphisms.   Note that sweepouts could be defined using a Lipschitz family of biLipschitz maps in place of diffeomorphisms.   The width is achieved by a minimal surface which may be unstable.  For example the width of a round sphere ${\mathbb{S}}^3$ is $4\pi$ and it is achieved by any equatorial sphere.  The equatorial sphere is unstable because it can be globally deformed to a surface of smaller area.  See Figure~\ref{fig-sphere}.  Marques and Neves proved the following rigidity theorem for spheres in \cite{MN-Duke}:

\begin{thm}{\bf Scalar Sphere Rigidity [Marques-Neves]}
Suppose $M^3$ is a sphere with a metric that has $\Scal \ge 6$ and $\Width(M^3)=4\pi$ then $M^3$ is isometric to a round sphere, ${\mathbb{S}}^3$.  
\end{thm}

This theorem is proven using Hamilton's Ricci Flow.  Hamilton proved that under Ricci flow scalar curvature remains positive and evolves in a beautiful way \cite{Hamilton-JDG82}.  
Colding and Minicozzi studied the evolution of Width
under Ricci flow to prove finite time extinction of the Ricci flow.  See for example their expository article \cite{CoMin-width}.    
Marques and Neves then applied the Ricci flow in combination with minmax theory to prove their rigidity theorem stated above in \cite{MN-Duke}.

\subsection{Schoen-Yau and Zero Mass}

\begin{figure}[h]
\begin{center}
\includegraphics[width=0.95\textwidth]{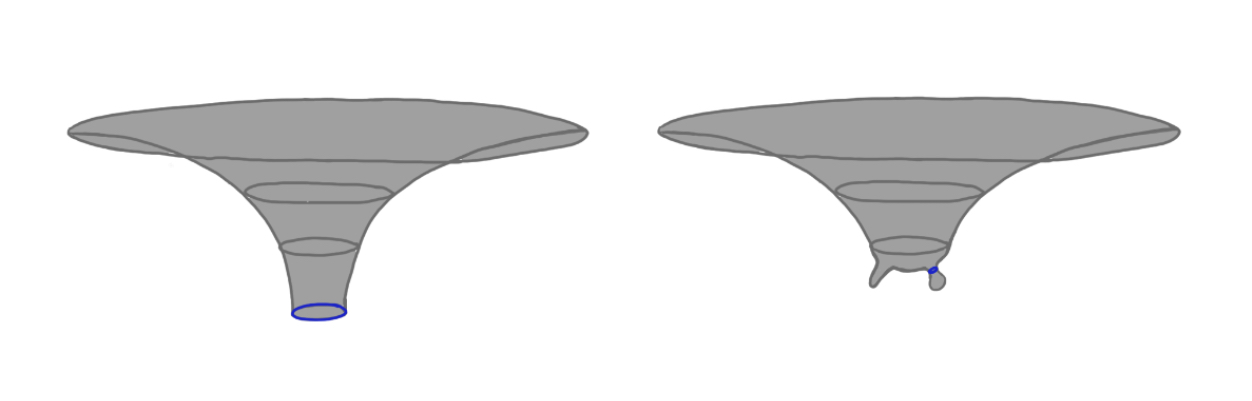} 
\caption{On the left is Riemannian Schwarzschild space which has $\Scal=0$  everywhere and a minimal boundary which is the horizon of a black hole with $m_{ADM}(M)= \sqrt{\Area(\partial M)/(16\pi)}$. 
The manifold on the right has $\Scal\ge 0$ and the same ADM mass due to its asymptotic behavior at infinity.  It has a well and a minimal surface in the neck of a bubble.}\label{fig-ADM-mass}
\end{center}
\end{figure}

Schoen and Yau's original Positive Mass Theorem states that a complete noncompact Riemannian manifold, $M^3$, which is asymptotically flat that has nonnegative scalar curvature has positive ADM mass unless the manifold is isometric to Euclidean space \cite{Schoen-Yau-PMT}.     See Figure~\ref{fig-ADM-mass}.    Although we are considering compact manifolds here, the techniques applied to prove this theorem and quasilocal versions of this theorem, may be applied in the compact setting as well.

While the ADM mass was originally defined in \cite{ADM}
using a high level of regularity, Huisken introduced an isoperimetric mass
in \cite{Huisken-Isoper} which requires only notions of volume and area.  The ADM mass has recently been shown by Miao
using volume estimates of Fan-Shi-Tam 
in \cite{FST09} to be equal to Huisken's isoperimetric mass:
\be \label{mISO}
m_{ADM}(M^3)=m_{ISO}(M^3)=\limsup_{r\to \infty}( \Omega_r)
\ee
where $\Omega_r$ are increasingly large round regions and 
\be
m_{ISO}(\Omega) = \frac{2}{\Area(\partial \Omega)} \left( \vol(\Omega) - \frac{\Area(\partial \Omega)^{3/2}}{6\sqrt{\pi}}\right).
\ee
See also work of Jauregui-Lee in \cite{Jauregui-Lee-isoper} and Chodosh-Eichmair-Shi-Yu
in \cite{CESY}.

Here we will state the rigidity part of the positive mass theorem in dimension three as follows:

\begin{thm} {\bf Zero Mass Rigidity [Schoen-Yau]}
If $M^3$ is asymptotically flat with $\Scal \ge 0$ everywhere and $m_{ADM}(M^3)=0$
then $M^3$ is isometric to Euclidean Space.
\end{thm}

Witten proved the Positive Mass Theorem for Spin manifolds in \cite{Witten-PMT}.  Lohkamp has 
another approach outlined in \cite{Lohkamp-PMT}.   It should be noted that Lohkamp introduced a method of reducing the proof of the Positive Mass Theorem to the Torus Rigidity Theorem.   This is explained in more detail in Schoen-Yau's paper
on the Positive Mass Theorem in dimensions greater than 8 \cite{Schoen-Yau-high}.  
There is a Ricci flow proof by Yu Li in \cite{Li-RicciflowPMT}. Most recently there is
a new proof in three dimensions using harmonic maps to circles by Bray-Kazaras-Khuri-Stern \cite{BKKS-PMT}.

The Zero Mass Rigidity Theorem has been proven for smooth manifolds with corners by Miao \cite{Miao-corners} and by Shi-Tam in \cite{Shi-Tam-JDG02}.    More recent papers with even lower regularity have been written by 
McFeron-Sz\'{e}kelyhidi \cite{MS-2012}, Shi-Tam in \cite{Shi-Tam-singular}, Lee-LeFloch \cite{Lee-LeFloch},
and Jiang-Shen-Jang \cite{JSZ21}.  

Huisken-Ilmanen proved the Penrose Inequality:
\be\label{Penrose}
m_{ADM}(M^3) \ge \sqrt{\Area(\partial M)/(16\pi)\,}
\ee
if $\partial M$ is a connected minimal surface and there are no closed interior minimal surfaces in $M^3$
\cite{Huisken-Ilmanen-01}.  They proved the
Penrose Rigidity Theorem which states that if equality is achieved in (\ref{Penrose}) then $M^3$ is the
Riemannian Schwarzschild manifold depicted in Figure~\ref{fig-ADM-mass}.   The Penrose Inequality has been proven when $\partial M^3$ has more than one connected component by Bray using a completely different approach \cite{Bray-Penrose}.   It should be noted that the
Penrose Inequality fails if there are closed interior minimal surfaces surrounding the boundary of $M^3$.  In General Relativity closed minimal surfaces are apparent horizons of black holes, and it is well understood that anything can happen beyond the horizon of a black hole.   This intuition can prove helpful when studying compact manifolds with 
$\Scal\ge 0$ as well.

There are quasilocal notions of mass defined for regions with boundaries that are mean convex  with  positive Gauss curvature.  In addition to the isoperimetric mass mentioned above, there is the {\em Hawking mass}, which 
can take negative values even on regions in Euclidean space.   Christoudoulou-Yau have proven it is nonnegative on round regions in \cite{CY-remarks}.  Sun has proven rigidity of the boundaries of round regions with zero Hawking mass in \cite{Sun-PJM18} and Shi-Sun-Tian-Wei have proven rigidity of the interiors in \cite{SSTW-rigidity}.  Huisken-Ilmanen have proven the Hawking mass is nonnegative on regions defined using inverse mean curvature flow.  This flow is not defined for compact manifolds \cite{Huisken-Ilmanen-01}, however it can be applied to regions in compact manifolds which have extensions to complete noncompact manifolds with no closed interior minimal surfaces.  
The Bartnik mass of a region is defined as the infimum of the ADM mass over all such extensions
\cite{Bartnik-mass}.   

Finally there is the {\em Brown-York mass} which has been proven to be equal to
\be
m_{BY}(\partial \Omega)= \frac{1}{8\pi} \left(A_0'(0)-A'(0)\right)
\ee
where $A(r)$ is the area of the boundary of the tubular neighborhood,
\be
A(r) = \Area(\partial (T_r(\Omega)) )
\ee
and 
\be
A_0(r) =  \Area(\partial (T_r(\Omega_0)) )
\ee
where $\Omega_0\subset {\mathbb{E}}^3$ is a convex region in Euclidean space determined by the existence of a Riemannian isometric embedding from $\partial \Omega$ into Euclidean space that is unique up to isometry by a result of Nirenberg \cite{Nirenberg-embedding}
and Pogorelov \cite{Pogorelov-embedding}.   If the metric is only $C^1$ the isometric embedding may not be unique.  
To define the Brown-York mass in a lower regularity setting: we could perhaps choose among convex embeddings, an embedding with the maximum value for $A'_0(0)$ .

Shi and Tam have proven a rigidity theorem for the Brown-York mass in \cite{Shi-Tam-JDG02}. 
Imitating the statement of Gromov's Prism Rigidity, we state their theorem as follows:

\begin{thm} {\bf Zero Local Mass Rigidity  [Shi-Tam]}  
If $M^3$ has $\Scal \ge 0$, then any region $\Omega \subset M^3$ satisfies local mass rigidity.  That is:
if $\partial \Omega$ is mean convex and has positive Gauss curvature and $m_{BY}(\partial \Omega)\le 0$ 
then $\Omega$ is isometric to a region in Euclidean space and consequently $m_{BY}(\partial \Omega)= 0$.
\end{thm}

Shi and Tam proved this theorem  by extending the region $\Omega$ to a manifold $N$ containing $\Omega$
that is asymptotically flat.  They prove the extension has $m_{ADM}(N)=0$  and thus by the positive mass theorem with corners, the extension, $N$, is isometric to Euclidean space.  Thus the region $\Omega\subset N$  is a flat region in Euclidean space.   The theorem has been extended by Eichmair, Miao, and Wang to more general shaped regions using a different proof in \cite{Eichmair-Miao-Wang}.

\begin{rmrk}
A smooth manifold $M^3$ whose small balls satisfy local mass rigidity has nonnegative scalar curvature.   This is
a consequence of the fact that
\be \label{lim-BY}
\lim_{r\to 0} 12 m_{BY}(\partial B(p,r))/r^3= \Scal(p)
\ee
proven by Fan, Shi, and Tam 
for a manifold with a $C^2$ metric tensor in \cite{FST09}.  
\end{rmrk}  

\begin{quest} If a $C^0$ manifold satisfies local mass rigidity for small balls, is the limit of the ratio of volumes in (\ref{Scal-vol-lim}) nonnegative? 
\end{quest}

\begin{quest} 
Is there a direct relationship between local mass rigidity and prism rigidity?   
\end{quest}

\begin{quest}
Can we define a natural quasilocal mass for prisms which depends only on areas, mean curvature, and dihedral angles? 
\end{quest}

The above two questions were discussed at IAS and are being investigated by various mathematicians.   See in particular recent work of Miao in  \cite{Miao-cubes} 
and by Miao-Piubello in \cite{Miao-Piubello}.

\begin{rmrk} \label{rmrk-Miao-vol} 
All the quasilocal masses are defined depending on the mean curvature of the boundary of $\Omega$, which is why they are written as $\mass(\partial \Omega)$ rather than $\mass(\Omega)$.  Suppose we have a surface $\Sigma$ with prescribed mean curvature.  It is natural to ask whether it has a filling $\Omega$ with nonnegative scalar curvature whose boundary 
with the induced metric is isometric to $\Sigma$.  One consequence of Shi-Tam's Rigidity Theorem is that such a filling only exists if the resulting Brown-York mass is positive using the integral of the prescribed mean curvature in the place of $A'(0)$.  Miao and Tam note in \cite{Miao-Tam-CVPDE09} that it is unknown what the smallest volume of such a possible filling might be even when $\partial \Omega$ is isometric to the boundary of a Euclidean ball.   They observe in this case that the Euclidean filling is definitely not the one of smallest volume because the volume functional
is unstable at the Euclidean filling.  See 
also their joint work with Jauregui in \cite{JMT13} and 
the more recent work of Yuguang Shi, Wenlong Wang, Guodong Wei, and Jintian Zhu in \cite{SWWZ20}. 
This question is discussed further by Gromov in \cite{Gromov-Natural}. 
\end{rmrk}

The rigidity of asymptotically hyperbolic manifolds with $\Scal \ge -6$ were first considered by Min-Oo in \cite{Min-Oo-89}.
The mass of such manifolds is studied in Wang in \cite{Wang-JDG01} and Chrusciel-Herzlich \cite{CH03}.  
The rigidity has been proven in work of Andersson-Dahl \cite{Andersson-Dahl-98},
of Andersson-Cai-Galloway \cite{ACG-Rigidity}, and of Huang-Jang-Martin \cite{HJM-Rigidity}.

Other rigidity theorems for complete noncompact manifolds with scalar curvature bounds are the Hyperbolic Scalar Cylinder Rigidity Theorem of Nunes in \cite{N-JGA-13} and Moraru in \cite{M-JGA16} and the Rigidity of Regions between Minimal surfaces of Carlotto-Chodosh-Eichmair \cite{CCE-Effective}.   There is also the Cylindrical Splitting Theorem of Chodosh-Eichmair-Moraru in \cite{CEM-Splitting} and the Warped Scalar Splitting Theorem of Galloway-Jang in \cite{Galloway-Jang-20}.  See also the Vacuum Static Rigidity Theorem of
Qing-Yuan \cite{Qing-Yuan-vacuum}.  Finally there is the
rigidity of manifolds with bad ends proven by Lesourd-Unger-Yau in \cite{LUY-20}.

This concludes our very limited survey of the geometric properties and geometric rigidity theorems for three dimensional Riemannian manifolds with scalar curvature bounds.   These results are the ones we will be discussing below in relation to convergence, but keep in mind there are many other beautiful results concerning manifolds with scalar curvature bounds which will be surveyed in other chapters of this volume.   It should be noted that we only consider lower bounds on scalar curvature throughout because Lohkamp has proven in \cite{Lohkamp-Ann95} that any Riemannian manifold can be achieved as the $C^0$ limit of sequences of manifolds with scalar curvature bounded uniformly from above.

\section{Geometric Notions of Convergence} \label{Sect-SWIF}

Here we will present the following notions of convergence reviewing their geometric properties: 
\begin{itemize}
\item Gromov-Lipschitz ($Lip$) convergence 
\item Intrinsic flat ($\mathcal{F}$) convergence 
\item Volume preserving intrinsic flat ($\mathcal{VF}$) convergence 
\item Measure convergence ($m$) 
\item Volume above distance below (VADB) convergence 
\end{itemize}   
We will describe how these notions are related to one another and how they relate to less geometric notions of convergence defined by the convergence of metric tensors.    We may view this section as a very brief review of the
key ideas and a survey of the existing literature.   

\subsection{Lipschitz Convergence}
We begin with Lipschitz convergence of compact metric spaces as defined by Gromov in \cite{Gromov-text}, not to be confused with the Lipschitz convergence of metric tensors.   
First recall that the Lipschitz constant of a function $F: (X_0,d_0)\to (X_1, d_1)$ is defined
\be
\Lip(F) = \sup_{p\neq q} \frac{ d_1(F(p),F(q))}{d_0(p,q)}
\ee
Note that
\be
d_1(F(p),F(q))\, \le \, \Lip(F) \, d_0(p,q) 
\ee
immediately implies that the rectifiable length of a curve,
\be \label{rect-length}
L_d(C[a,b]) = \sup\left\{ \sum_{i=1}^N d(C(t_i), C(t_{i-1}))\, : a=t_0<\cdots <t_N=b,\,\, N \in {\mathbb{N}} \,\right\},
\ee
satisfies
\be \label{Lip-length}
L_{d_1}((F\circ C)[a,b]) \le \Lip(F) L_{d_0}(C([a,b])
\ee
and similarly the Hausdorff measures of sets satisfy
\be \label{Lip-measure}
\mathcal{H}^m_{d_1}(F(U)) \,\le\,  (\Lip(F))^m \,\mathcal{H}^m_{d_0}(U) .
\ee

We say that a sequence of metric spaces $(M_j,d_j)$ converges to $(M_\infty, d_\infty)$ in the
Gromov-Lipschitz sense \cite{Gromov-text},
\be\label{defn-Gromov-Lip}
(M_j,d_j) \LIPto (M_\infty, d_\infty)
\ee
if there exists biLipschitz maps, $\Psi_j: M_\infty \to M_j$, such that
\be
\max\left\{\Lip(\Psi_j), \Lip(\Psi^{-1}_j)\right\} \le 1+\epsilon_j \to 1.
\ee

Given a fixed manifold with a sequence of metric tensors $g_j \to g_0$ in the $C^0$ sense then 
\be \label{C^0-to-Lip}
(M, d_j) \LIPto (M, d_\infty)
\ee
where $d_j$ are defined as in (\ref{d-from-g}).   In fact the Lipschitz constant of the identity map 
from $(M, d_j)$ to $(M, d_\infty)$ is 
bounded above by $\max_{v\in TM} g_0(v,v)/g_j(v,v)$ and vice versa.

Note that under Lipschitz convergence, it is easy to see that 
\be
 \Diam(M_j) \to \Diam(M_\infty)
 \textrm{ and }
 \Vol(M_j) \to \Vol(M_\infty).
 \ee
So we say Lipschitz convergence is {\em diameter preserving} and {\em volume preserving}.   

If $M_j$ has a boundary, $\partial M_j$, then we could restrict the distances on $M$ to distances on the boundary.   With the restricted distances, we get Lipschitz convergence of the boundaries.  That is Lipschitz convergence is {\em boundary preserving}:
\be
(\partial M_j, d_j) \to (\partial M_\infty, d_\infty)
\ee
and thus also {\em preserves the area of the boundary}:
\be
\Area(\partial M_j) \to \Area(M_\infty).
\ee
Within the paper there will be further discussion of restricted vs induced distances on the boundaries.  One advantage of using restricted distances is that the boundary can have many components and all are considered part of the same disconnected metric space.   For more information about the convergence of manifolds with boundary see
the survey by Perales \cite{Perales-survey16}.

\begin{figure}[h]
\begin{center}
\includegraphics[width=0.95\textwidth]{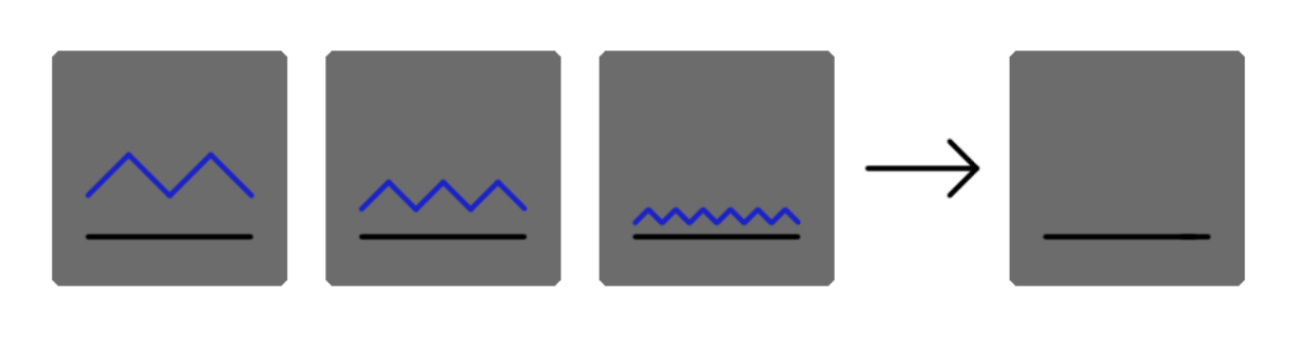} 
\caption{The sequence of jagged curves, \textcolor{blue}{$C_j$}, converge to a straight line segment,
$C_\infty$, but $\lim_{j\to \infty} L_j(C_j) > L_\infty(C_\infty)$.}\label{fig-drop-length}
\end{center}
\end{figure}

Lipschitz convergence {\em preserves lengths} in the following sense:
any rectifiable curve in $C: [0,1]\to M_\infty$ has a sequence of rectifiable curves $\Psi_j \circ C:[0,1] \to M_j$ which "converge" to this limit curve and whose lengths converge:
\be
L_j(\Psi_j \circ C) \to L_\infty(C)
\ee
as a consequence of (\ref{Lip-length}).
Conversely, if we start with
curves $C_j:[0,1] \to M_j$ with a uniform upper bound on their lengths, the Arzela-Ascoli Theorem implies that they have a subsequence which converges to rectifiable curves in $C_\infty: [0,1]\to M_\infty$ in the sense that $\Psi_{j_k}^{-1}\circ C_{j_k}: [0,1] \to M_\infty$
converge uniformly to $C_\infty: [0,1]\to M_\infty$.  We even have
\be
\liminf_{k\to \infty} L_{j_k}(\Psi_{j_k} \circ C) \ge L_\infty(C)
\ee
but the lengths may drop in the limit and the limit curve might even be a single point.   This can be seen by considering
a constant sequence $M={\mathbb{T}}^2$ with $g_j =g_0$ the standard flat metric, and $\Psi_j$ the identitiy map and taking a sequence of $C_j$ as in Figure~\ref{fig-drop-length}.  

\begin{figure}[h]
\begin{center}
\includegraphics[width=0.95\textwidth]{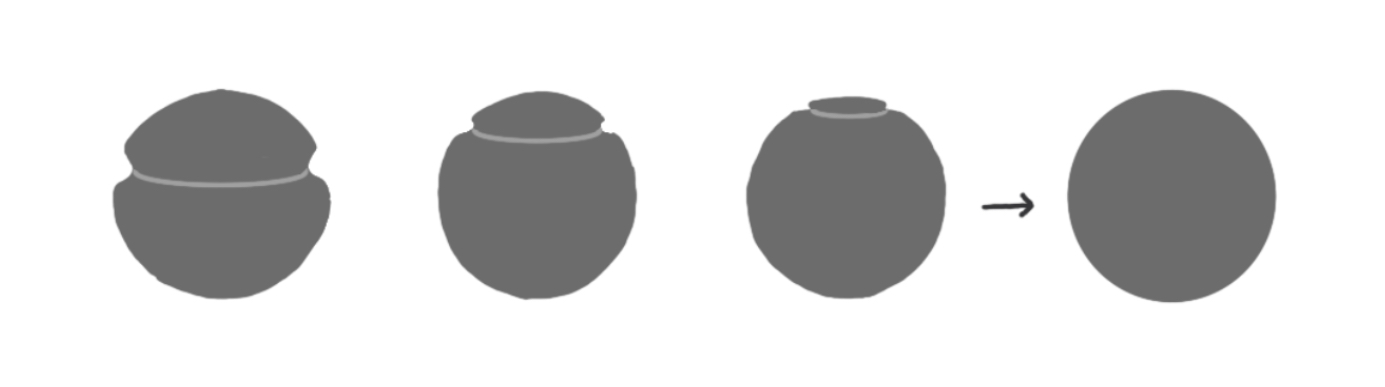} 
\caption{This sequence of warped two spheres are Lipschitz converging to a standard sphere.  They have stable minimal geodesics that disappear in the limit due to shortening length.   These spheres can also be viewed as a sequence of warped three spheres with stable minimal surfaces that disappear in the limit.}\label{fig-latitudes-to-zero}
\end{center}
\end{figure}

Let us now consider rectifiable surfaces, $\Sigma \subset M_\infty^3$.   As above we have surfaces
$\Psi_j(\Sigma) \subset M_j$ and the Hausdorff measures converge:
\be
\Area(\Psi_j(\Sigma)) \to \Area(\Sigma).
\ee
In this sense Lipschitz convergence is {\em area preserving}.  Conversely, if we have closed surfaces $\Sigma_j\subset M_j$, it is possible for them to disappear in the limit, even if they are stable minimal surfaces.  For example if $M_j= [0,\pi] \times_{f_j} {\mathbb{S}^2}$ and $f_j(t)$ converge to $\sin(t)$
in the $C^0$ sense then $M_j \LIPto {\mathbb S}^3$.  Choose a sequence of such $f_j$ such that
$f_j'(t_j)=0$ and $f_j"(t_j)>0$ we have $\Sigma_j=f_j^{-1}(t_j)$ is a stable minimal sphere.   If we choose $t_j\to 0$
we have $\Area(\Sigma_j)\to 0$ and the $\Sigma_j$ disappear in the limit.  See Figure~\ref{fig-latitudes-to-zero} .   Even if the $\Sigma_j$ don't disappear, they might converge to a surface which is no longer stable as seen by taking $t_j \to \pi/2$ or a surface which is no longer minimal as seen by taking $t_j \to \pi/4$.   This will be discussed more within the paper, as will a discussion of the limits of mean convex surfaces and constant mean curvature surfaces.

If we assume that we have oriented rectifiable surfaces $\Sigma_j\subset M_j$ such that
\be
\Area(\Sigma_j) \le A \textrm{ and } \Length(\partial \Sigma_j) \le L
\ee
then
\be
\Area(\Psi_j^{-1}(\Sigma_j)) \le A' \textrm{ and } \Length(\Psi_j^{-1}(\partial \Sigma_j)) \le L'.
\ee
Thus we may apply the Federer-Flemming Compactness Theorem \cite{FF}, to conclude that a subsequence converges in the flat (FF) sense as integral currents to a limit 
\be
\Lbrack\Psi_j^{-1}(\Sigma_j)\Rbrack \FFto T_\infty
\ee
where $T_\infty$ is an integral current, which either corresponds to a weighted oriented countably $\mathcal{H}^2$ rectifiable surface or is the $0$ current.  

Let us explain this more clearly.  A subset $S$ of a metric space, $(M,d)$,  is said to be {\em countably $\mathcal{H}^m$ rectifiable} and oriented if it has a preferred collection of countably many disjoint biLipschitz charts, $\varphi_i: A_i \to U$ defined on Borel sets, $A_i \subset {\mathbb E}^m$, such that
\be
\mathcal{H}^m\left(\, S \, \setminus \, \bigcup_{i=1}^\infty \varphi_i(A_i) \, \right) = 0.
\ee
We can integrate $m$ dimensional forms, $\omega$, over such a set to define a {\em current}:
\be
\Lbrack S \Rbrack (\omega)= \int_S \omega = \sum_{i=1}^\infty \int_{A_i} \varphi_i^* \omega.
\ee
If we add integer weights $a_i \subset {\mathbb N}$, we obtain what is called an {\em integer rectifiable current}:
\be
T(\omega)= \sum_{i=1}^\infty a_i \Lbrack \varphi_i(A_i) \Rbrack \textrm{ where } T(\omega)=\sum_{i=1}^\infty a_i \int_{A_i} \varphi_i^* \omega.
\ee
When $M$ is a Riemannian manifold, $T$ has a {\em mass} which is its {\em weighted volume}:
\be
\mass(T) = \sum_{i=1}^\infty |a_i| \mathcal{H}^m(\varphi_i(A_i)).
\ee
The {\em boundary} of a current is defined abstractly by Stokes Theorem
\be
(\partial T)(\omega)= T(d\omega)
\ee
so that if $\Sigma$ is a smooth surface with boundary we have 
$\partial \Lbrack \Sigma \Rbrack= \Lbrack \partial \Sigma\Rbrack$
the two notions of boundary agree with the correct orientation.

All these ideas have been extended to metric spaces by Ambrosio-Kirchheim \cite{AK} which we will describe within.
The formula for mass of a rectifiable current in a metric space is
more complicated than just a weighted volume.  Neither notion of mass has anything to do with ADM mass or
quasilocal mass.
 
\begin{figure}[h]
\begin{center}
\includegraphics[width=0.95\textwidth]{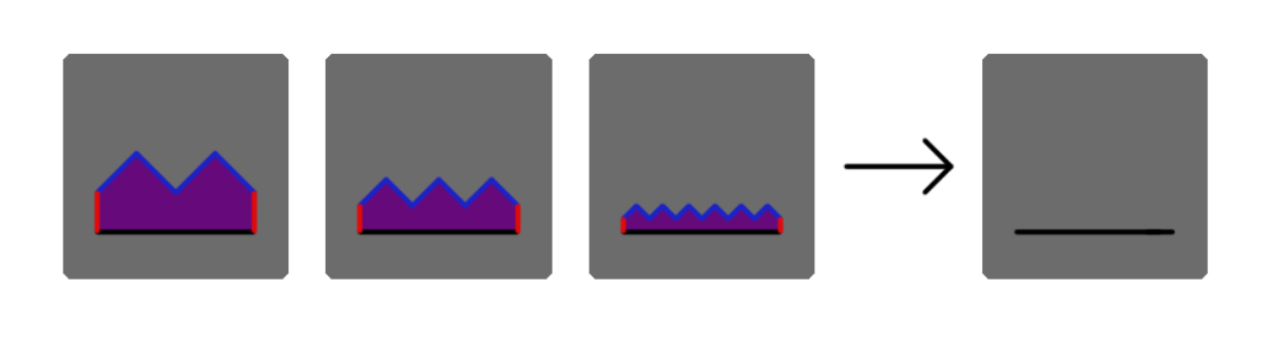} 
\caption{The sequence of jagged curves are represented by currents, \textcolor{blue}{$T_j =\Lbrack C_j \Rbrack$}, which
F-converge to the line segment,
$T_\infty=\Lbrack C_\infty \Rbrack $.  This can be seen by taking \textcolor{violet}{$B_j =\Lbrack \Omega_j \Rbrack$} where $\Omega_j$ is the region between the connected curves and  \textcolor{red}{$A_j$} to be the vertical segments on either side
oriented so that $T_\infty-T_j-A_j = \partial B_j$.}\label{fig-drop-length-F}
\end{center}
\end{figure}

An {\em integral current} is an integer rectifiable current that has a boundary that is also integer rectifiable.  We say a sequence of integral currents $T_j$ converges in the flat (F) sense as integral currents to $T_\infty$ if
\be
d_F^M(T_j, T_\infty) \to 0
\ee
where flat distance between currents is defined as in Federer-Flemming \cite{FF} as
\be
d_F^M(T_j, T_\infty)= \inf\{\mass(A) + \mass(B)\,:\,
T_j - T_\infty= A + \partial B \, \}
\ee
where $A$ is also an $m$ dimensional integral current and $B$ is an $m+1$ dimensional integral current
whose boundary is built from $T_j$, $T_\infty$, and $A$ with appropriate orientations defined by the signs above.
See Figure~\ref{fig-drop-length-F}.   Note that if $T_j$ is a sequence such that $\mass(T_j)\to 0$ then $T_j$ converges to the $0$ current.   

It is easy to see from the definition that $T_j \FFto T_\infty$ implies the boundaries converge:
$
\partial T_j \FFto \partial T_\infty.
$
Furthermore the masses are lower semicontinuous:
\be
\liminf_{j\to \infty}\mass(T_j) \ge \mass(T_\infty); 
\ee
so the same is also true for the boundaries:
\be
\liminf_{j\to \infty}\mass(\partial T_j) \ge \mass(\partial T_\infty). 
\ee
Another beautiful property of flat convergence that will be described further within is the {\em convergence of slices}.  If we have a Lipschitz function $h: M \to \mathbb{R}$ and a rectifiable surface $\Sigma\subset M$
then we can define slices $<h,\Lbrack \Sigma \Rbrack, t>$ which correspond to the sets $h^{-1}(t) \cap \Sigma$ for almost 
every $t$.   If $\Lbrack \Sigma_j \Rbrack \FFto \Lbrack \Sigma_\infty\Rbrack $ then for almost every $t \in {\mathbb R}$,
a subsequence of the slices converge 
\be
<h,\Lbrack \Sigma_{j_k}\Rbrack , t>\,\, \FFto \,\, <h,\Lbrack \Sigma_\infty \Rbrack, t>.
\ee
Federer-Flemming's flat convergence has been a fundamental part of the study of Geometric Measure Theory since the 1950's and we recommend Morgan's textbook for an intuitive introduction \cite{Morgan-GMT}.  We recommend Ambroisio-Kirchheim's {\em Currents on Metric Spaces} \cite{AK} for a rigorous and elegant extension of the theory to metric spaces that works beautifully as an introduction to the theory on Euclidean Space as well.   

\subsection{Intrinsic Flat Convergence}

As mentioned in the introduction, Lipschitz convergence of Riemannian manifolds is too strong a notion of convergence to consider for sequences of $M_j^3$ with nonnegative scalar curvature.   We will still consider it as an important test case for checking which properties of scalar curvature might be conserved when taking a limit of such manifolds, but it is not a weak enough notion of convergence to apply to consider more general sequences of manifolds with nonnegative scalar curvature.   

Sormani and Wenger weakened the notion of Lipschitz convergence of oriented Riemannian manifolds using the idea of flat convergence of integral currents in \cite{SW-JDG}.   A sequence of compact oriented Riemannian manifolds with boundary converges in the intrinsic flat ($\mathcal{F}$) sense
\be\label{SWIF}
M^n_j \Fto M^n_\infty \iff d_{\mathcal{F}}\left(M^n_j, M^n_\infty\right) \to 0
\ee
where
\be
d_{\mathcal{F}}\left(M_j, M_\infty\right)
=\inf \,\,d_F^Z\left( \Lbrack \Psi_j(M_j) \Rbrack, \Lbrack \Psi_\infty(M_\infty)\Rbrack \right),
\ee
where the infimum is taken over all complete metric spaces, $(Z,d_Z)$, and all distance preserving maps
\be
\Psi_j: M_j \to Z \textrm{ such that } d_Z\left(\Psi_j(p),\Psi_j(q)\right)=d_j(p,q) \quad \forall p,q \in M_j.
\ee
They use Ambrosio-Kirchheim theory to rigorously define the integral currents, $\Lbrack\Psi(M)\Rbrack$,
on the metric space $Z$.   This $\mathcal{F}$ convergence is also referred to as Sormani-Wenger Intrinsic Flat (SWIF) convergence in some later papers.   
See Figure~\ref{VF-Example} for a sequence of manifolds which converge in the intrinsic flat sense to a sphere and note that there is no requirement that the topology of the manifolds may change in the limit.   

Sormani and Wenger proved that  
\be
d_{\mathcal{F}}(M_j, M_\infty)=0
\ee
iff there is an orientation preserving isometry between $M_j$ and $M_\infty$.
 They also proved that if $M_j \Fto M_\infty$ then there exists a common complete metric space $Z$, and
distance preserving maps $\Psi_j: M_j \to Z$ such that 
\be\label{common-Z}
\Lbrack \Psi_j(M_j) \Rbrack \FFto \Lbrack \Psi_\infty(M_\infty)\Rbrack \textrm{ as integral currents in } Z.
\ee
Thus intrinsic flat convergence has all the same properties as flat convergence:
\be
M_j \Fto M_\infty \implies \partial M_j \Fto \partial M_\infty
\ee
and
\be
\liminf_{j\to \infty} \Vol(M_j) \ge \Vol(M_\infty) \textrm{ and } \liminf_{j \to \infty} \Vol(\partial M_j) \ge \Vol(\partial M_\infty).
\ee
In \cite{S-ArzAsc} Sormani has also proven that
\be
\liminf_{j\to \infty} \Diam(M_j) \ge \Diam(M_\infty).
\ee

Sormani and Wenger defined a larger class of oriented rectifiable weighted metric spaces, called
{\em integral current spaces}, which include the $0$ space, and defined the $\mathcal{F}$ convergence of such spaces
in \cite{SW-JDG}.   Oriented Alexandrov spaces are integral
current spaces as seen in work of Jaramillo, Perales, Rajan, Searle, Siffert, and Mitsuishi  \cite{JPRSS}\cite{Mit-coin}\cite{Mit-orient}.   Integral Current Spaces
need not be connected nor have geodesics, and they even include the $0$ space \cite{SW-JDG}.   
A sequence of Riemannian manifolds $\mathcal{F}$ converges to the $0$ space if they are collapsing,
\be
M_j \Fto 0 \textrm{ if } \vol(M_j) \to 0,
\ee
or if they are cancelling as depicted in an example in \cite{SW-JDG}.      Two beautiful properties about integral current spaces are that they have boundaries
(which are also integral current spaces) satisfying
\be
M_j\Fto M_\infty \implies \partial M_j \Fto M_\infty
\ee
and at almost every point they have a unique tangent cone which is a normed space by \cite{AK}.
For integral current spaces it is the Ambrosio-Kirchheim mass of the space which is semicontinuous rather than the volume:
\be
M_j \Fto M_\infty \implies \liminf_{j\to \infty} \mass(M_j) \ge \mass(M_\infty)
\ee
This mass is more than just a weighted volume as it also has an area factor that depends on the structure of the tangent spaces.  Sormani and Wenger prove that if the intrinsic flat distance between two integral current spaces is $0$ then there is a current preserving isometry between them which preserves the orientation, the integer-valued weights, and the area factors \cite{SW-JDG}.

In \cite{Wenger-compactness}, Wenger proved that if a sequence of Riemannian manifolds has
\be \label{Wenger-compactness-eqn}
\Diam(M_j) \le D \textrm{ and } \Vol(M_j) \le V \textrm{ and } \Area(\partial M_j) \le A
\ee
then a subsequence converges in the $\mathcal{F}$ sense to an integral current space which is possibly the $0$ space (c.f. \cite{SW-JDG}).   It should be noted that $\mathcal{F}$ limits do not necessarily agree with Gromov-Hausdorff limits.
Even when a Gromov-Hausdorff limit exists, $M_j \GHto M_{GH}$, the $\mathcal{F}$ limit is only known to
be a subset $M_\infty \subset M_{GH}$.  Since $M_\infty$ is always the same dimension as $M_j$, we know $M_\infty$ is the $0$ space when $M_{GH}$ is lower dimensional \cite{SW-JDG}.   

In \cite{PS-properties}, Portegies and Sormani have proven that the filling volume is continuous:
\be
M_j \Fto M_\infty \implies \FillVol(M_j) \to \FillVol(M_\infty)
\ee
where the filling volume is defined using the Ambrosio-Kirchheim mass:
\be
\FillVol(M):= \inf\{ \mass(N) \, | \, (\partial N, d_N)=(M,d_M)\}.
\ee
Here the infimum is taken over all
integral current spaces $N$ whose boundary with the restricted distance from $N$ has a current preserving
isometry to M.    Gromov's Filling Volume defined in \cite{Gromov-filling} might have a different value.   The relationship between these two filling volumes is not well studied yet, but Gromov's theorems seem likely to carry over and so far have been easily adapted and proven as needed to apply to this Ambrosio-Kirchheim mass filling volume.   Note that there is a wealth of literature further developing and applying Ambrosio-Kirchheim Theory that can be consulted when working with their notion of mass.     There is also a sliced filling volume defined by Portegies and Sormani in \cite{PS-properties} which behaves well under $\mathcal{F}$ convergence. 

\begin{figure}[h]
\begin{center}
\includegraphics[width=0.95\textwidth]{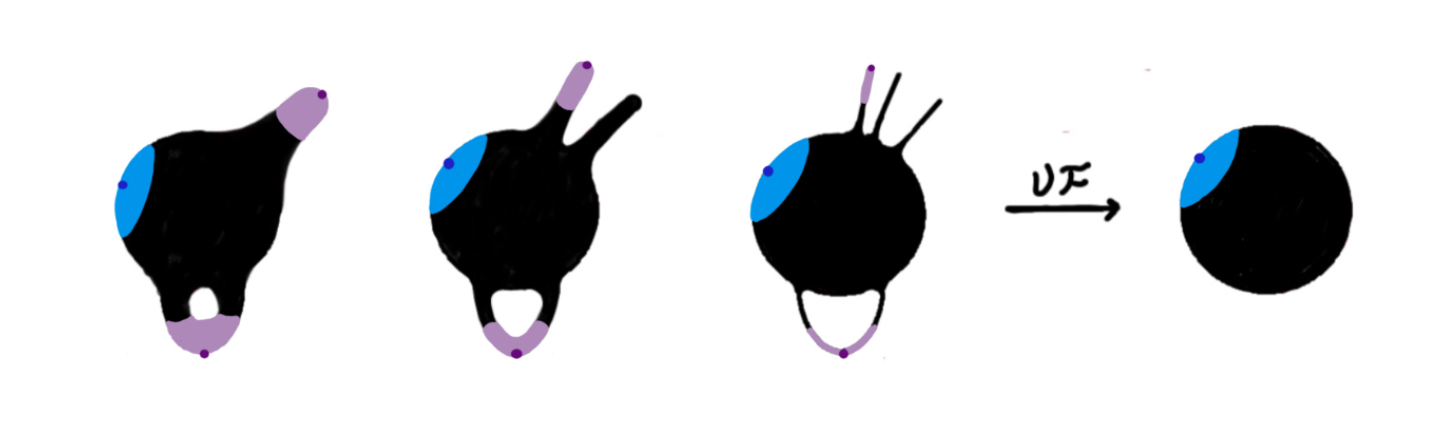} 
\caption{If $M_j \Fto M_\infty$ and \textcolor{blue}{$p_j \to p_\infty$} then for almost every $r>0$ a subsequence \textcolor{blue}{$B(p_{j_k},r)\Fto B(p_\infty,r)$.}
If  \textcolor{violet}{$q_j$} lie on collapsing regions then \textcolor{violet}{$B(q_j,r)\Fto 0$} and \textcolor{violet}{$q_j$} disappear in the limit. }\label{fig-balls-converge}
\end{center}
\end{figure}

We can localize the above theorems and apply them to balls in our converging spaces.
In \cite{S-ArzAsc}, Sormani proved that if a sequence of points $p_j\subset M_j$ converge to $p_\infty\subset M_\infty$ in the sense that $\Psi_j(p_j) \to \Psi_\infty(p_\infty)$ for $\Psi_j$ satisfying (\ref{common-Z}) then
for almost every $r>0$ there is a subsequence such that
\be
B(p_{j_k},r) \Fto B(p_\infty,r). 
\ee
Thus we have semicontinuity of the volumes of balls and areas of their boundaries.  See Figure~\ref{fig-balls-converge}.

Combining this with the work of Portegies and Sormani in \cite{PS-properties}, we have continuity of the filling volumes of balls and spheres about converging points.  In fact, we can prevent the disappearance of points $p_j$ by proving that 
\be\label{fill-balls}
\exists R>0 \,\, \textrm{ s.t. for a.e. } r \in (0,R) \,\,\, \exists C_r>0 \textrm{ s.t. }\FillVol\left(\partial B(q_j,r)\right)\ge C_r.
\ee
Applications of this technique will be described further below.

Sormani also proves in \cite{S-ArzAsc} that any $p\in M_\infty$ has a sequence $p_j \to p_\infty$ in the above sense.  In fact there are regions $U_j \subset M_j$ which GH converge to $M_\infty$.
Sormani applies this to prove two Arzela-Ascoli type theorems.    If $M_j$ and $X$ are compact and
\be
M_j \Fto M_\infty\neq 0 \textrm{ and } F_j: M_j \to X \textrm{ have } \Lip(F_j) \le K
\ee
then a subsequence converges to a Lipschitz function $F_\infty: M_\infty \to X$ where
$F_\infty(p_\infty)=\lim_{j\to \infty} F_j(p_j)$ for a sequence $p_j \to p_\infty$.   If sequences of compact spaces
with compact limits
 \be
M_j \Fto M_\infty\neq 0 \textrm{ and } N_j \Fto N_\infty \neq 0
\ee
and $F_j: M_j \to N_j$ is an isometry on balls of uniform radius $R>0$, then a subsequence converges to a local isometry.
Both of these theorems require an assumption that the spaces have a uniform upper bound on their volume and on the areas of their boundaries.

\begin{figure}[h]
\begin{center}
\includegraphics[width=0.95\textwidth]{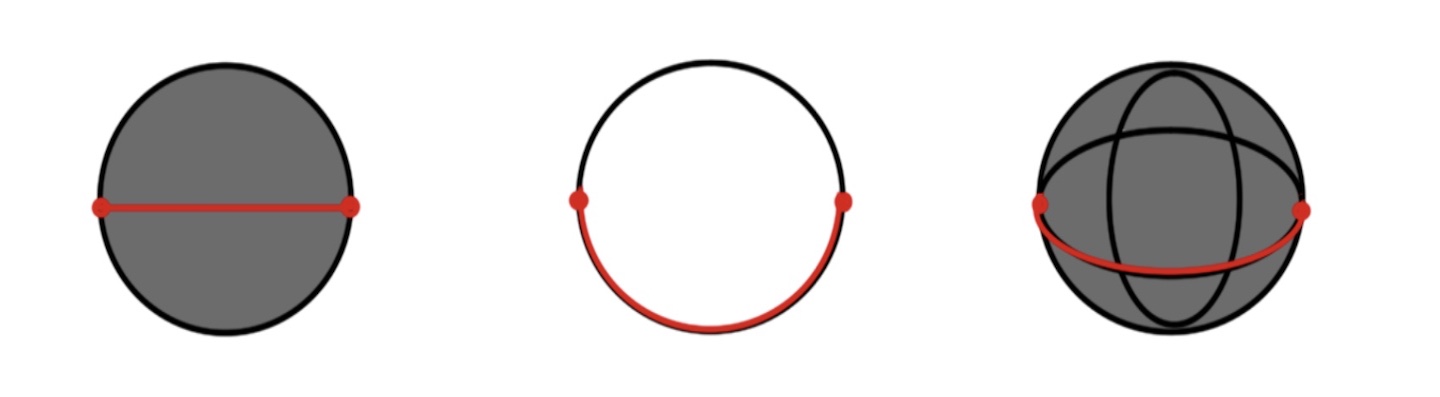} 
\caption{The map to the left  from ${\mathbb S}^1$ into a closed Euclidean disk is not distance preserving
but the map to the right into a longitude of ${\mathbb S}^2$ is distance preserving.}\label{fig-dist-pres}
\end{center}
\end{figure}

To estimate the intrinsic flat distance between $M_j$ and $M_\infty$ and prove convergence, we need only construct a sequence of metric spaces $Z_j$, find distance preserving maps $\Psi_j: M_j \to Z_j$ and
$\Psi_j': M_\infty \to Z_j$ and then measure the volume between their images.   It is important to note that the distance preserving maps cannot just be Riemannian isometries.   Riemannian isometries are only length preserving, do not guarantee that we have
\be
d_Z(\varphi_j(p),\varphi_j(q)) = d_j(p,q) \qquad \forall p,q \in M_j
\ee
as it is possible that the points in $Z$ could be joined by a short path in $Z$ that does not lie in
$\varphi_j(M_j)$.   For example the map $\varphi: M={\mathbb{S}}^1 \to Z={\mathbb E}^2$ that maps the circle to a unit circle in the Euclidean plane is not distance preserving but $\varphi: M={\mathbb{S}}^1 \to {\mathbb{S}}^2$
that maps the circle to the equator is distance preserving.   See Figure~\ref{fig-dist-pres}.

Lee and Sormani have proven a theorem estimating the intrinsic flat distances between manifolds when there are Riemannian isometric embeddings of those manifolds into a common Riemannian manifold, $N$, which are not distance preserving \cite{LeeSormani1}.   Their theorem, which was applied to study $M_j$ with $\Scal \ge 0$, requires estimating an embedding constant 
\be
C_M= \sup_{p,q\in M} d_M(p,q) - d_N(\varphi(p), \varphi(q))
\ee
and is proven by explicitly gluing together a metric space 
\be
Z=(M_i\times [0,h_i]) \disjointunion_{\varphi_i(M_i)} N
\ee
where $h_i$ depends on $C_{M_i}$ \cite{LeeSormani1}.

Using a different construction,  Lakzian and Sormani proved in  \cite{Lakzian-Sormani} that
\be
(M, g_j) \LIPto (M, g_\infty) \implies (M, g_j) \Fto (M, g_\infty)
\ee
by taking $Z_j= M \times [0,h_j]$ endowed with a metric tensor created using a combination of the metric tensors $g_j$ and $g_\infty$ so that 
\be
\varphi_j: (M, d_j) \to M\times \{h_j\} \subset Z_j \textrm{ and }
\varphi'_j: (M, d_\infty) \to M\times \{0\} \subset Z_j 
\ee
are both distance preserving.  They set $B_j=Z_j$ and $A_j=0$, so that
\be
d_{\mathcal{F}}( (M, g_j), (M,g_\infty)) \le \mass(B_j) +0 =\Vol(Z_j) \to 0.
\ee
In the same paper, Lakzian and Sormani proved a frequently applied theorem in which only open sets 
$U_j\subset M_j$ and $U_\infty\subset M_\infty$
are required to be bi-Lipschitz close.   They prove $M_j \VFto M_\infty$ as long as
there exists $\phi_j: U_j \to U_\infty$ whose biLipsschitz constants converge to 1, constants
\be
\lambda_j=\max\{\,|d_j(p,q)-d_\infty(\phi_j(p),\phi_j(q))|\, :  p,q \in U_j\} \to 0,
\ee
\be
\Vol(M_j \setminus U_j) \to 0,
\ee
and there are uniform upper bounds $\Vol(M_j)$, $\Area(\partial U_j)$, and $\Diam(M_j)$.
In their paper, this theorem was applied to prove sequences of manifolds converging smoothly away from a singular set have a variety of possible $\mathcal{F}$ limits.   Applications to study $M_j$ with $\Scal(M_j) \ge 0$
will be surveyed further within.

\subsection{Volume Preserving Intrinsic Flat and Measure Convergence}

We say that a sequence of compact oriented Riemannian manifolds without boundary converges in the 
{\em volume preserving intrinsic flat} sense 
\be
M_j \VFto M_\infty \iff M_j \Fto M_\infty \textrm{ and } \Vol(M_j) \to \Vol(M_\infty).
\ee
When $M_j$ are only  integral current spaces, we require mass convergence rather than volume converging.
The additional requirement that the volumes or masses converge is unexpected powerful and has important
consequences.

Portegies has proven that 
\be
M_j \VFto M_\infty \qquad \implies M_j \mto M_\infty
\ee
in \cite{Portegies-evalue}.   We say that a sequence of Riemannian
manifolds {\em converges in measure} to a limit space $M_j \mto M_\infty$ if and only if 
there exists a common complete metric space $Z$, and
distance preserving maps $\Psi_j: M_j \to Z$ such that 
\be
||\Lbrack\Psi_{j}(M_j)]]|| \to ||\Lbrack \Psi_{\infty}(M_\infty) \Rbrack || \textrm{ weakly as measures in } Z.
\ee
That is, for any $\mathcal{H}^m$ measurable set  $A\subset Z$ we have
\be
||\Lbrack \Psi_{j}(M_j)\Rbrack||(A)=\Vol_j \left(\Psi^{-1}_j(A)\right) \to ||\Lbrack\Psi_{\infty}(M_\infty)\Rbrack||(A)=\Vol_\infty \left(\Psi^{-1}_\infty(A)\right)
\ee
as long as $||\Lbrack\Psi_{\infty}(M_\infty)\Rbrack ||(\partial A)=0$ as in Portmanteau's Theorem.   
For more information about convergence in measure see Lott-Villani and Sturm's papers \cite{Lott-Villani} \cite{Sturm} and the
more recent work of Ambrosio, Gigli, Mondino, and Savare \cite{Ambrosio-Gigli-Savare} \cite{GMS}.  Note that convergence in measure does not imply
$\mathcal{VF}$ convergence as it does not control orientations nor boundaries.   

Portegies  then applied his theorem and techniques of Fukaya \cite{Fukaya-evalue} to prove that the {\em Laplace spectrum} of the $M_j$ is semicontinuous under $\mathcal{VF}$ convergence:  
\be \label{Portegies-evalue}
M_j \VFto M_\infty \qquad \implies \qquad \limsup_{j\to \infty} \lambda_k(M_j) \le \lambda_k(M_\infty)
\ee
where $\lambda_k(M)$ is the $k^{th}$ eigenvalue of the Laplacian counting multiplicity when $M$ is a Riemannian
manifold \cite{Portegies-evalue}.   If $M$ has boundary then $\lambda_k(M)$ is a Neumann eigenvalue. Portegies also describes what happens when the limit space is only an integral current space.  
Note that Fukaya's work in  \cite{Fukaya-evalue} concerned measured Gromov-Hausdorff convergence in which the sequence converges in measure and in the GH sense, while the $\mathcal{VF}$ sequences studied by Portegies need not even have GH limits.   

\begin{figure}[h]
\begin{center}
\includegraphics[width=0.5\textwidth]{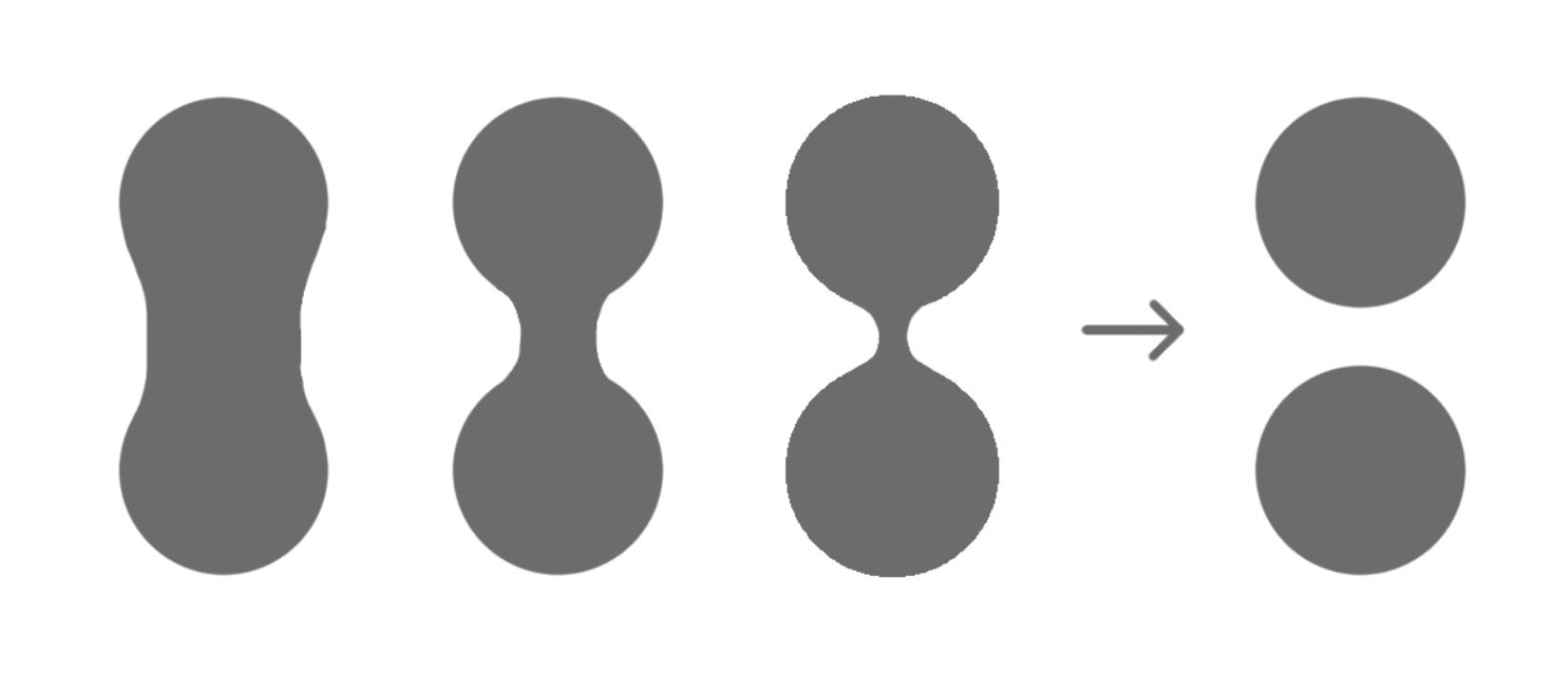} 
\caption{This sequence of warped three spheres $(M_j,g_j)$ with $\Scal_j \ge 0$, $\mathcal{VF}$ converge to a pair of disjoint spheres $M_\infty={\mathbb{S}}^3 \disjointunion  {\mathbb{S}}^3$ with a distance $d_\infty$ that is restricted from the GH limit which is a length metric on two spheres joined by a line segment.  }\label{fig-two-spheres}
\end{center}
\end{figure}

In general GH and $\mathcal{VF}$ limits do not agree.    However, Matveev-Portegies and Honda have proven in \cite{Matveev-Portegies-17} and \cite{Honda-Ricci-17} that when 
\be \label{Ricci}
\Ricci(M^n_j) \ge -(n-1)H \textrm{ and } \vol(M^n_j) \ge V_0>0 \textrm{ and } \partial M^n_j =0
\ee
then the GH, measured, and $\mathcal{VF}$ limits all agree.  With appropriate bounds on the boundary,
Perales has proven they agree when $\partial M^n_j$ is nonempty but has uniform bounds in
\cite{Perales-JTA20}.
  If we replace Ricci curvature in (\ref{Ricci}) by $\Scal(M_j)>0$ then one can create counter examples.  The simplest would be a sequence of manifolds formed
by taking a pair of three dimensional spheres joined by an increasingly thin tunnel as in Figure~\ref{fig-two-spheres}.
Such tunnel constructions were first described in the work of Schoen-Yau \cite{Schoen-Yau-structure} and Gromov-Lawson \cite{Gromov-Lawson-classification}.   The GH limit of such a sequence is a pair of spheres with a line segment between them, while the $\mathcal{VF}$ limit of the sequence is just the pair of spheres with the distance restricted from
the GH limit.   For a detailed construction of a tunnel of arbitrary length and width see the work of Dodziuk
\cite{Dodziuk-tunnels}.

Jauregui and Lee have done further work on the properties of 
$\mathcal{VF}$ convergence and its relationship with scalar curvature in \cite{Jauregui-Lee-VF}.   They've proven that if a sequence of points $p_j\subset M_j$ converge to $p_\infty\subset M_\infty$ in the sense that $\Psi_j(p_j) \to \Psi_\infty(p_\infty)$ for $\Psi_j$ satisfying (\ref{common-Z}) then
for almost every $r>0$
\be
B(p_j,r) \VFto B(p_\infty,r) 
\ee
and so we have continuity of the volumes of balls.  Note it is still possible for balls to disappear but only if $\Vol(B(q_j,r))\to 0$ as depicted in Figure~\ref{fig-balls-converge}.  Jauregui and Lee have also proven continuity of the areas of their boundaries and also of the areas of level sets of other Lipschitz functions  \cite{Jauregui-Lee-VF}.  They applied their 
work to prove the semicontinuity of ADM mass under intrinsic flat convergence.

\subsection{VADB Convergence}

Allen, Bryden, Perales, and Sormani have completed a series of papers relating $GH$, $\mathcal{F}$, and $\mathcal{VF}$ convergence to various notions of convergence of metric tensors on a fixed Riemannian manifold
\cite{AS-relating}\cite{AS-contrasting}\cite{Allen-Bryden-19}\cite{VADB} \cite{Allen-Perales-VADB}.   This includes examples demonstrating that $L^p$ convergence of metric tensors does not necessarily imply $\mathcal{F}$ convergence of the manifolds \cite{AS-relating}.  Nevertheless Allen, Perales, and Sormani were able to prove that if we have $L^p$ convergence of
metric tensors for $p \ge n/2$
and $g_j \ge (1-(1/j)) g_\infty$ then we have $(M^n,d_j) \VFto (M^n, d_\infty)$ \cite{VADB}.    

We say a sequence of compact oriented Riemannian manifolds without boundary has {\em volume above distance below} (VADB) {\em convergence}: $M_j \VADBto M_\infty$ if
\be \label{VADB}
\Vol(M_j) \to \Vol(M_\infty) \textrm{ and } \Diam(M_j) \le D  
\ee
and there exists $C^1$ diffeomorphisms $\Psi_j: M_\infty \to M_j$ with
\be \label{dist-below}
d_j(\Psi_j(p),\Psi_j(q)) \ge d_\infty(p,q)\quad \forall (p,q) \in M_\infty\times M_\infty.
\ee
Note that $M_j \VADBto M_\infty$ when $M_j=(M,g_j)$ with $g_j \ge g_\infty$ satisfying (\ref{VADB}).
We can even rescale $M_j$ to study manifolds with $g_j \ge g_0 -1/j$ satisfying (\ref{VADB}).

Allen, Perales, and Sormani have proven in \cite{VADB} that 
\be
M_j \VADBto M_\infty \quad \implies M_j \VFto M_\infty.
\ee
The first step towards proving this theorem appeared in \cite{AS-contrasting} as Theorem 2.6, where Allen and Sormani showed
$M_j \VADBto M_\infty $ implies
\be \label{VADB-ptwise}
 \lim_{j\to \infty} d_j\left(\Psi_j(p),\Psi_j(q)\right) \ge d_\infty(p,q)
\textrm{ pointwise a.e. } (p,q) \in M_\infty\times M_\infty.
\ee
Allen-Perales-Sormani  then prove there exists $\delta_j \to 0$ and $V_j \to 0$ and domains
$W_j \subset M_\infty$ such that
\be
\Vol_j\left(M_j \setminus \Psi_j(W_j)\right) \le V_j \textrm{ and }
\ee
\be
\delta_j \ge d_j\left(\Psi_j(p),\Psi_j(q)\right) - d_\infty(p,q) \ge 0 \qquad \forall \, p,q\in W_j
\ee
Finally they estimate the intrinsic flat distance
\be
d_{\mathcal{F}}(M_j,M_\infty)\le 2 V_j +  V \sqrt{2 \delta_j D + \delta_j^2}
\ee
by constructing a metric space as follows:
\be
Z_j= \left(M_\infty \times\{0\}\right) \disjointunion \left(W_j \times [0,h_j]\right) \disjointunion \left(M_j \times \{h_j\} \right).
\ee
and embedding both $M_j$ and $M_\infty$ into the top and bottom of $Z_j$ with distance preserving maps.  
  
Allen and Perales have proven similar theorems for VADB convergence of manifolds with boundary \cite{Allen-Perales-VADB}.   They
have a variety of different hypotheses we may assume on the boundary.  The simplest version is to assume the 
interiors of the manifolds are convex and the areas of the boundaries are uniformly bounded above.   
In each of their theorems, they prove their hypotheses combined with (\ref{VADB}) and (\ref{dist-below})
imply $M_j \VFto M_\infty$.   

Applications of VADB convergence to prove $\mathcal{VF}$ convergence of sequences of manifolds
with lower bounds on their scalar curvature  appear in work of Cabrera Pacheco, Ketterer, and Perales \cite{CKP-torus-graph}, work of Huang, Lee, and Perales \cite{HLP}, and work of Allen \cite{Allen-conf-torus}.   These will be discussed in more detail below.

\section{MinA Scalar Compactness [Gromov-Sormani]}\label{Compact1}

\begin{figure}[h]
\begin{center}
\includegraphics[width=0.5\textwidth]{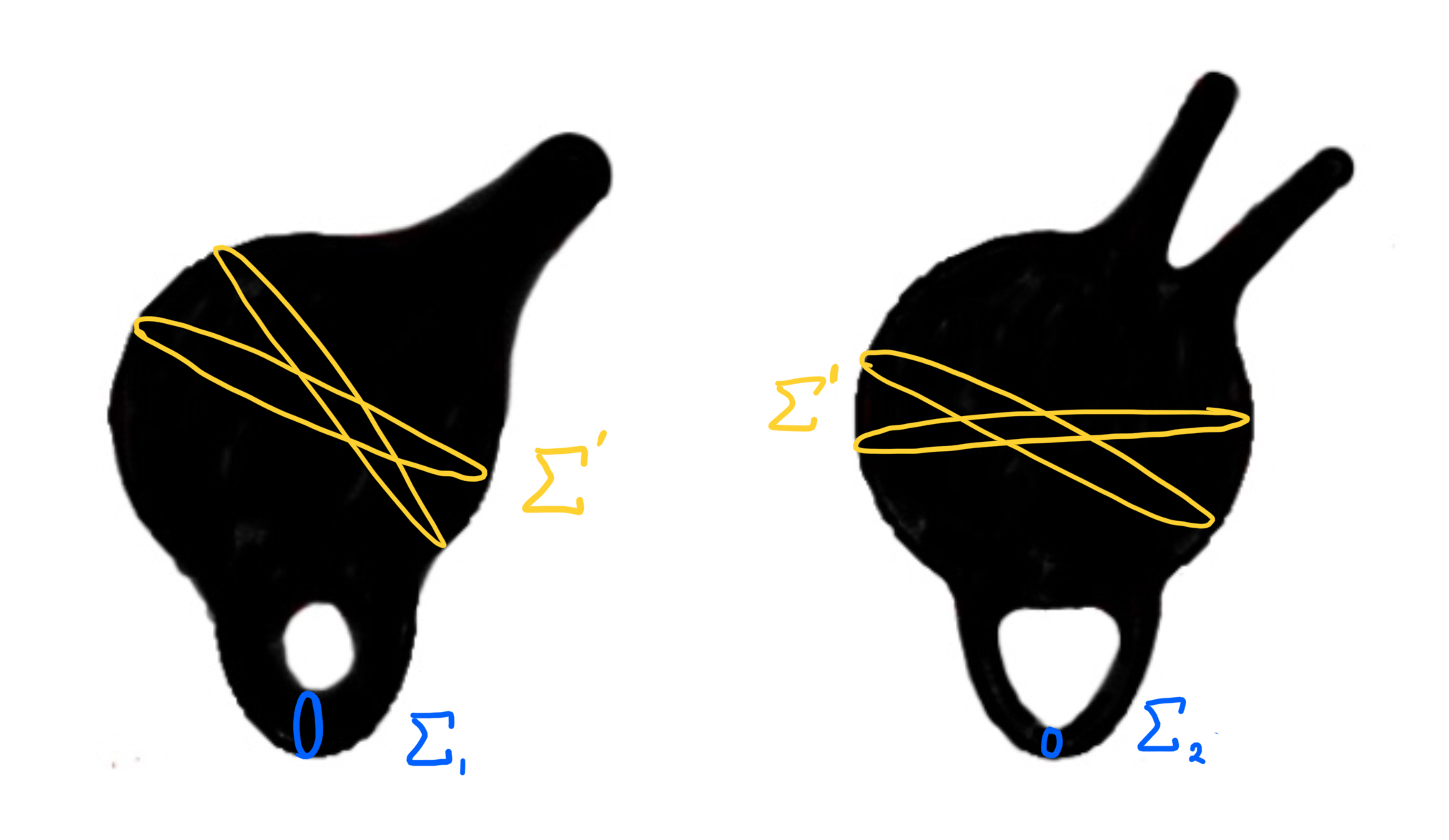} 
\caption{The minimal spheres $\Sigma_i$ are stable and $\Sigma'$ are unstable.   }\label{minsurfaces}
\end{center}
\end{figure}

The first Scalar Compactness Conjecture was suggested by Gromov in \cite{Gromov-Plateau} and further refined by Sormani in \cite{Sormani-scalar} building upon work with Basilio and Dodziuk in \cite{BDS-sewing} and work of Park-Tian-Wang in \cite{Park-Tian-Wang-18}. 
First recall the definition:
\be \label{MinA}
\MinA(M^3) =\min\{ \Area(\Sigma)\, | \, \Sigma \textrm{ is a closed minimal surface in } M^3 \}.
\ee
Recall that a minimal surface is locally area minimizing.  See Figure~\ref{minsurfaces} for images of stable minimal surfaces which are globally area minimizing and unstable minimal surfaces which are not.   Schoen-Yau proved that stable minimal surfaces in manifolds with nonnegative scalar curvature are spheres or tori \cite{Schoen-Yau-minimal}.

\begin{conj} {\bf{MinA Scalar Compactness Conjecture [Gromov-Sormani]}} \label{Conj-MinA}
Suppose $M_j^3$ are oriented three dimensional Riemannian manifolds with 
\be \label{ScalComp1}
\Scal_j \ge 0 \textrm{ and } \MinA(M_j^3) \ge \alpha >0
\ee
satisfying
\be \label{ScalComp2}
\Diam(M_j) \le D \textrm{ and } \Vol(M_j) \le V \textrm{ and }  \partial M_j =\emptyset
\ee
then a subsequence $M_{j_k}^3 \VFto M_\infty$ where $M_\infty$ is a 
three dimensional rectifiable geodesic metric space with Euclidean tangent cones
satisfying generalized notions of nonnegative scalar curvature like the Prism Rigidity Property
or the Local Mass Rigidity Property.   
\end{conj}

Note by Wenger's Compactness Theorem as stated in \cite{SW-JDG} and
the hypotheses in (\ref{ScalComp2}) imply that
a subsequence that converges in the $\mathcal{F}$ sense.   The first challenge then is to prove that the limit
is not the $0$ space.    Without the $\MinA$ condition the sequence $\vol(M_j^3)$ could collapse to the $0$ space.  For example,
a sequence of increasingly thin flat tori
\be
M_j^3={\mathbb{S}}_{1/j}^1 \times {\mathbb{S}}^1\times {\mathbb{S}}^1
\ee
which has 
\be
\MinA(M_j^3)= \Area({\mathbb{S}}^1_{1/j} \times {\mathbb{S}}^1)=(2\pi/j)(2\pi)=4 \pi^2/j \to 0.
\ee
Note that the definition of $\MinA$ must include unstable minimal surfaces as can be seen by considering
\be
M_j^3={\mathbb{S}}^1_{1/j} \times {\mathbb{S}}^2
\ee
which has minimal surfaces diffeomorphic to ${\mathbb{S}}^1_{1/j} \times {\mathbb{S}}^1$
of area $2\pi/j \times 2\pi \to 0$.   

Even if the sequence has a uniform lower bound on volume, without the $\MinA$ condition, the
$M_j$ can $\mathcal{VF}$-converge to a pair of disconnected spheres and thus fail to be a geodesic metric space.  
See Figure~\ref{fig-two-spheres}.   An example by Basilio, Kazaras, and Sormani demonstrates that without the $\MinA$ condition the limit space might fail to have any geodesics at all.   Indeed they produce a sequence of $M_j^3$
with $\Scal>0$ satisfying (\ref{ScalComp2}) that $\mathcal{VF}$ converge to a three sphere with the restricted metric from Euclidean space \cite{BKS-no-geod}.   See Remark~\ref{rmrk-geods} for a discussion as to how we  might to prove the limit
$M_\infty$ of Conjecture~\ref{Conj-MinA} is geodesic.

Conjecture~\ref{Conj-MinA} has been proven in the warped product setting by Park-Tian-Wang in \cite{Park-Tian-Wang-18}.  They assume the $M_j$
are three dimensional spheres and that the metric tensor has the form $dt^2 + f_j^2(t) g_0$ where $f(0)=0$
and $f_j(L_j)=0$ and $f_j(t)>0$ between.   Using the various hypotheses they were able to
prove that a subsequence of the $f_j$ converges in the $H^1_{loc}$ sense to $f_\infty \ge 0$ which
is positive on an interval $(a,b)$.   They then prove 
\be
M_j \VFto M_\infty \textrm{ where } M_\infty=({\mathbb S}^3, dt^2 + f_\infty^2(t) g_0)
\ee
They prove that the scalar curvature is nonnegative in a generalized sense defined using the fact that this is
an $H^1_{loc}$ metric tensor.  LeFloch and Sormani had previously defined nonnegative scalar curvature in the same way when investigating a related compactness problem involving Hawking mass \cite{LeFloch-Sormani}.

A more simply stated definition of scalar curvature in a manifold of lower regularity is to say that
\be \label{Scal-vol-liminf}
\Scal(p)=\liminf_{r\to 0} 6(n+2) \left(\frac{V^n(r) - \Vol(B(p,r))}{r^2 V^n(r)} \right).
\ee
In \cite{BDS-sewing}, Basilio, Dodziuk, and Sormani construct sequences of closed Riemannian
manifolds, $M_j^3$,
satisfying all the hypotheses of Conjecture~\ref{Conj-MinA} except the $\MinA$ condition, which $\mathcal{VF}$
and $mGH$ converge to a space, $M_\infty$, which has a point where this limit is $-\infty$.   The example is constructed by taking a curve in the sphere and {\em sewing along the curve} with a sequence of precisely spaced increasingly tiny increasingly dense tunnels so that the limit $M_\infty$ is a sphere with the curve identified to a point.    In \cite{BS-seq}, Basilio and Sormani prove an entire region can be sewn to a point by locating the endpoints of the tunnels in a precisely described way within the region.
  
\begin{quest}\label{Quest-scal-lim}
If we have a sequence $M_j^3$ satisfying all the hypotheses of Conjecture~\ref{MinA} and we assume the
sequence converges in the $C^0$, Lipschitz, $m$, GH, $\mathcal{VF}$, $VADB$, or $H^1_{loc}$  
sense, is 
\be \label{Scal-vol-liminf}
\liminf_{r\to 0} 6(n+2) \left(\frac{V^n(r) - \Vol(B(p,r))}{r^2 V^n(r)} \right) \ge 0
\ee
at every $p$ or perhaps at almost every $p$?   It is known for $C^0$ convergence if we  assume the limit space is
smooth by work of Gromov and Bamler mentioned above.   
\end{quest}

\begin{rmrk} Gromov has suggested that we should require that the limit space in our compactness theorem satisfy the Scalar Prism Rigidity Theorem for all prisms that can be drawn in the space.   In \cite{Gromov-Dirac}, Gromov has checked that under $C^0$ convergence, the Scalar Prism Rigidity Property is conserved.     It would be interesting to check if a Riemannian
manifold, $M_\infty$, that is a VADB limit of a sequence of Riemannian manifolds, $M_j$,
satisfying the hypotheses of the above conjecture  satisfies the Scalar Prism Rigidity Property as well.
For more general limit spaces, we must ensure that dihedral angles and the notion of mean convexity is preserved under convergence.     
\end{rmrk}

\begin{rmrk} We might also ask whether the limit spaces satisfy the Local Mass Rigidity Property.    This property might be easier to define on general limit spaces than the Prism Rigidity Property because we would not need to define dihedral angles or find cubes.  However we would need to find a way to embed the boundaries of small balls in
Euclidean space.    It would be interesting to first check whether this property is conserved under Lip or even VADB convergence if we know the limit is a smooth Riemannian manifold.    
\end{rmrk}

\begin{rmrk}  
Bamler has proven in \cite{Bamler-C} that nonnegative scalar curvature is conserved under $C^0$ convergence of metric tensors applying a theorem of Simon concerning the properties of Ricci flows \cite{Simon02}.
using Ricci flow methods and 
\cite{Bamler-C}.   
His method might possibly be applied to verify if nonnegative scalar curvature is conserved under VADB convergence
of $C^2$ Riemannian manifolds to a $C^2$ Riemannian manifold.
\end{rmrk}

\begin{rmrk}
Burkhardt-Guim has applied Ricci flow starting at a $C^0$ manifold to define the scalar curvature
at a point in a $C^0$ manifold \cite{B-G-GAFA} .   It is unknown how her notion relates to the $H^1_{loc}$
definition of scalar curvature or the (\ref{Scal-vol-liminf}) 
definition but it is quite powerful.
She has proven this notion is conserved under $C^0$ convergence.   
How well does her notion behave under VADB convergence?   
\end{rmrk}

\begin{rmrk} \label{rmrk-min-appear}
A first step towards proving Conjecture~\ref{Conj-MinA} would be to prove by contradiction that a sequence satisfying the
hypotheses cannot $\mathcal{VF}$ converge to the $0$ space.   We could try to prove that when this happens, a minimal surface must appear.   In \cite{Schoen-Yau-existence}, Schoen and Yau prove that a minimal surface can be forced to appear within a region with a strong concentration of positive scalar curvature.   Andersson and Metzger have a theorem
providing upper bounds on the areas of minimal surfaces in \cite{AM-CMP09} that might be applied to contradict the
$\MinA$ hypotheses.
\end{rmrk}

\begin{rmrk}
This conjecture is three dimensional because we have no reason to believe minimal surfaces will be strong enough to prevent the collapse of higher dimensional manifolds.  There are important four dimensional examples of Man Chun Lee, Aaron Naber, and Robin Neumayer that would need to be better understood before stating a higher dimensional 
version of this conjecture.
In their paper they prove that assuming a uniform lower bound on entropy suffices to prove that there is a $d_p$-converging subsequence $M_{j_k}^m$ for $m\ge 3$.  It would be interesting to explore how their $d_p$ convergence relates to $\mathcal{F}$ convergence and how their entropy bound relates to $\MinA$ in dimension three.  However in dimension four, they provide a sequence of manifolds satisfying their theorem's hypotheses for which the $\mathcal{F}$ limit is the $0$ space \cite{Lee-Naber-Neumayer-dp}. 
\end{rmrk}

\begin{rmrk}
Li-Mantoulidis present intriguing examples in dimension 3 of Riemannian manifolds with positive scalar curvature away from skeleton singularities, which can be desingularized to smooth Riemannian manifolds with positive scalar curvature \cite{Li-Mantoulidis-skeleton}. They conjecture that their sequences converge in the pointed VF sense, but generally not in the Gromov-Hausdorff sense. They also conjecture such desingularization works for Riemannian manifolds with skeleton singularities in all dimensions (currently, this is proved for edge-cone singularities, a special class of skeleton singularities). 
\end{rmrk}

Let us consider the various consequences that would follow from Conjecture~\ref{Conj-MinA}.  Each of these conclusions leads to natural questions that might be examined on its own to find a counter example to the conjecture.   Proving one of the conclusions directly might have interesting applications that could be pursued even while the full conjecture remains open.   

The simplest to state is that there is a uniform lower bound on volume
for sequences of manifolds satisfying its hypotheses.   This follows because the conjecture concludes that the limit
space, $M_\infty$ is a nonzero space, so it has mass in the sense of Ambrosio-Kirchheim, and under
$\mathcal{F}$ convergence:
\be
\liminf_{j\to \infty} \Vol(M_j) \ge \mass(M_\infty) >0.
\ee
This conclusion does not even require $\mathcal{VF}$ convergence.   All we would need to prove is $\mathcal{F}$
convergence to a nonzero limit.
This leads to a natural question:

\begin{quest} \label{Quest-MinA-VolMin}
If we have a sequence $M_j^3$ satisfying all the hypotheses of Conjecture~\ref{MinA} 
can we find a lower bound on the volume of a manifold satisfying the
hypotheses of Conjecture~\ref{Conj-MinA} as a formula depending on
$\alpha$, $D$, and $V$?   Can we directly prove the lower bound is positive without
proving Conjecture~\ref{Conj-MinA} first? 
\end{quest}

Another consequence of $\mathcal{VF}$ convergence is the semiconvergence of the eigenvalues
which was proven by Portegies in \cite{Portegies-evalue}.   

\begin{quest} \label{Quest-MinA-eval-max}
If we have a sequence $M_j^3$ satisfying all the hypotheses of Conjecture~\ref{MinA} 
can we find an upper bound on the first eigenvalue of these manifolds as a formula depending only on
$\alpha$, $D$, and $V$?  Can we directly prove the upper bound is finite without
proving Conjecture~\ref{Conj-MinA}?  What can we say about the eigenmaps and heat kernels of the
manifolds?  
\end{quest}

\begin{quest} \label{Quest-MinA-fillvol}
If we have a sequence $M_j^3$ satisfying all the hypotheses of Conjecture~\ref{Conj-MinA}, then 
can we find a positive lower bound on the filling volume of these manifolds as a formula depending only on
$\alpha$, $D$, and $V$?     Note that very little is known about the lower bounds of filling
volumes.   It is not even known what the filling volume is for a round ${\mathbb{S}}^3$.   So it would perhaps
be easier to directly prove the lower bound is positive rather than to find a value for it.   
Finding such a lower bound would also be a possible first step towards proving Conjecture~\ref{Conj-MinA}, because a lower bound on the filling
volumes of a sequence guarantees that the sequence will not collapse to the $0$ space.  
\end{quest}

\begin{rmrk}
Note that Alexander Nabutovsky and Regina Rotman have theorems which state that when the Gromov Filling radii of three dimensional manifolds converge to $0$ then $MinA(M_j) \to 0$ in \cite{NR-GAFA-06}.  One would need to extend theorems of Gromov from \cite{Gromov-filling} and their results
to our filling volume defined using Ambrosio-Kirchheim mass.   
\end{rmrk}

\begin{rmrk}\label{rmrk-geods}
Let us consider how we might prove that the 
limit spaces, $M_\infty$, found in Conjecture~\ref{Conj-MinA} are geodesic metric spaces.  
If we wish to imitate Gromov's proof in \cite{Gromov-text} that GH limits are geodesic, we
start by taking an arbitrary pair of points  $p, p'\in M_\infty$ and prove they
have a midpoint $q\in M_\infty$, such that 
\be
d(p,p')=2d(p,q)=2d(q,p').
\ee
Gromov proved this by first finding $p_j, p'_j \in M_j$ converging to $p,p' \in M_{\infty}$ respectively
and we can do the same applying the work in \cite{S-ArzAsc}.   He then takes midpoints $q_j$
between $p_j$ and $p_j'$ and proves a subsequence of the $q_j$ converge to a midpoint $q$ between $p$ and $p'$
using the Bolzano-Weierstrass Theorem.   We can also find midpoints $q_j$
between $p_j$ and $p_j'$.   To obtain a limit for the midpoints, by \cite{PS-properties}, we need only to prove that
\be\label{fill-balls-2}
\exists R>0 \,\, \textrm{ s.t. for a.e. } r \in (0,R) \,\,\, \exists C_r>0 \textrm{ s.t. }\FillVol\left(\partial B(q_j,r)\right)\ge C_r.
\ee
Sequences of points, $q_j$, like those in Figure~\ref{fig-balls-converge}, can disappear in the limit and fail to satisfy (\ref{fill-balls}).
In that figure, the disappearing $q_j$ either lie in wells or in tunnels.  Those that are on the tips of wells are not midpoints of pairs of points, $p_j, p_j'$ that have limits $p,p'\in M_\infty$.   The disappearing $q_j$ that are in thin tunnels might be midpoints, but they are close to increasingly small minimal surfaces.   It would suffice to prove 
the following conjecture for $M_j^3 \VFto M_\infty$:
\be
\textrm{ if } B(p_j,R) \Fto B(p_\infty,R) \neq 0 \textrm{ and } B(p'_j,R) \Fto B(p'_\infty,R)\neq 0
\ee
and $q_j$ are midpoints of $p_j$ and $p_j'$ such that 
\be
\exists r>0 \textrm{ s.t. } \FillVol(B(q_j,r)) \to 0
\ee
then there exist minimal surfaces $\Sigma_j \subset M_j$
such that $\Area(\Sigma_j)\to 0$.  
\end{rmrk}

Conjecture~\ref{Conj-MinA}  itself is not likely to be proven for over a decade.  
This will become more clear as we discuss the various geometric stability theorems below, 
some of which will follow from a proof of this
conjecture with a sufficiently strong notion of generalized scalar curvature.  Solving any of the
questions above or proving anything suggested in one of the remarks is worthy of a paper.
It would also be worthwhile to prove the conjecture under additional hypotheses.   Some natural additional hypotheses 
are:
\begin{itemize}
\item $M_j$ are all conformal to a fixed manifold 
\item $M_j$ are all doubly warped product manifolds 
\item $M_j$ are all graphs over a fixed manifold 
\end{itemize}
The scalar curvature is well understood in these settings and we might even be able to prove VADB 
convergence which would then imply $\mathcal{VF}$ convergence.     Some of the geometric stability theorems have already been solved under these hypotheses as will be discussed below.  We recommend that anyone interested in proving one of these special cases to please contact the authors of the
corresponding geometric stability theorems as well as any one of us before proceeding to ensure that someone is not already working on the problem.   There are so many open questions in this area that there is no reason for people to work on the same problem with the same additional hypotheses unless they are working together.

\section{A Compactness Conjecture with Minimal Boundary}\label{Compact2}

\begin{figure}[h]
\begin{center}
\includegraphics[width=0.8\textwidth]{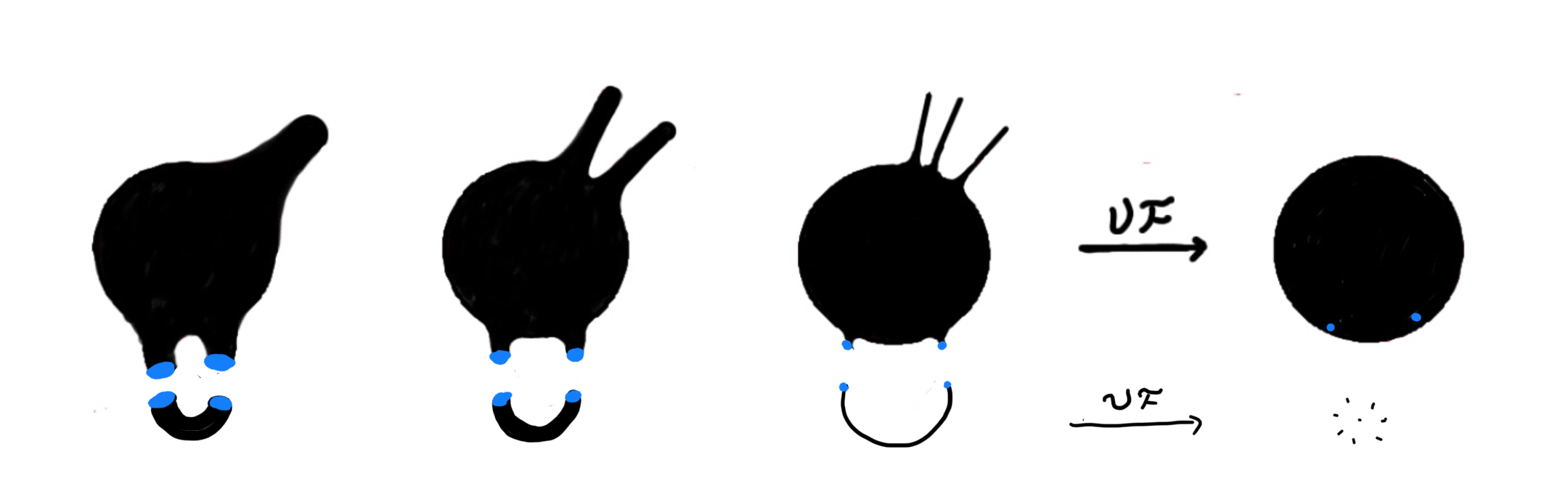} 
\caption{If we cut the manifolds from Figure~\ref{VF-Example} at stable minimal spheres, the regions on top converge to a sphere but the regions on the bottom collapse and disappear under $\mathcal{VF}$ convergence.}\label{cutmin}
\end{center}
\end{figure}

This next compactness theorem was formulated in conversations with Gromov at IAS and NYU in 2018 and has not been published before.  It addresses concerns that Gromov had that we might not be able to estimate $\MinA$ and so instead we might choose to cut up the sequence of manifolds along their stable minimal surfaces and study each region separately.   See Figure~\ref{cutmin}. The following conjecture could then be applied to the regions within the manifold.   As the hypotheses no longer prevents collapsing, it is possible the limit region might be the $0$ space.

\begin{conj}{\bf{NoMin Boundary Scalar Compactness [Gromov]}} \label{Conj-NoMin}
Suppose $M_j^3$ are oriented three dimensional Riemannian manifolds 
with boundaries such that
\be
\partial M_j \textrm{ are stable minimal surfaces }
\ee
such that
\be
\Scal_j\ge 0 \textrm{ and there are no closed interior minimal surfaces in $M_j$}
\ee
and furthermore
\be \label{W-Hyp}
\Diam(M_j) \le D \textrm{ and } \Vol(M_j) \le V \textrm{ and } \Area(\partial M_j) \le A 
\ee
then a subsequence $M_{j_k}^3 \VFto M_\infty$ where $M_\infty$ is either
the zero space or it is a
three dimensional rectifiable geodesic metric space with Euclidean tangent cones
satisfying generalized notions of nonnegative scalar curvature.
\end{conj}

Recall that Wenger's Compactness Theorem \cite{Wenger-compactness}  
the hypotheses in  (\ref{W-Hyp}) imply a subsequence converges to 
 $M_\infty$ that is either the zero space or a three dimensional rectifiable metric space called an integral current space (cf. \cite{SW-JDG}).   Note that by Ambrosio-Kirchheim theory, integral current spaces have tangent cones at almost every point that are normed vector spaces but they might not be Euclidean \cite{AK}.
So the fact that
a subsequence converges in the intrinsic flat sense is already known for both of these compactness conjectures.  The challenge is proving the volumes converge $\vol(M_j)\to \vol(M_\infty)$ to obtain $\mathcal{VF}$ convergence and
to prove the limit spaces are geodesic and satisfy the generalized notions of nonnegative scalar curvature.   Proving just one of these conclusions would be of interest even in a special case.

It is an important open question to find properties these limit spaces might satisfy that would allow us to prove the various stability conjectures we discuss below.  Note that while we do suggest a list of properties that might be useful, we do not list everything being explored in this direction by experts in scalar curvature.  In particular we do not discuss Spin approaches to defining scalar curvature here, although we do believe such an approach might prove fruitful. 

Although we have called these conjectures "compactness conjectures" they are in reality, precompactness conjectures
as the limit spaces are not Riemannian manifolds.  
This is not such a serious concern when discussing our first Conjecture~\ref{Conj-MinA}.  In that setting we would conjecture that if the limit happened to be a Riemannian manifold then it would satisfy the same hypotheses as the sequence approaching it.  Volume would be preserved by definition of $\mathcal{VF}$ convergence.  Diameter was proven to be lower semicontinuous in
\cite{S-ArzAsc}:
\be
M_j \VFto M_\infty \implies \liminf_{j\to \infty} \Diam(M_j) \ge \Diam(M_\infty).
\ee
Part of the conjecture is that the scalar curvature is nonnegative in some weak sense that implies nonnegative scalar curvature when $M_\infty$ is smooth.   Finally we would conjecture that
\be \label{limsup-MinA}
M_j \VFto M_\infty \implies \limsup_{j\to \infty} \MinA(M_j) \le \MinA(M_\infty)
\ee
under the hypotheses of Conjecture~\ref{Conj-MinA}, although it is unlikely that (\ref{limsup-MinA}) without the
scalar curvature bound in light of examples appearing in work of Sinaei-Sormani \cite{Sinaei-Sormani}.

\begin{figure}[h]
\begin{center}
\includegraphics[width=0.8\textwidth]{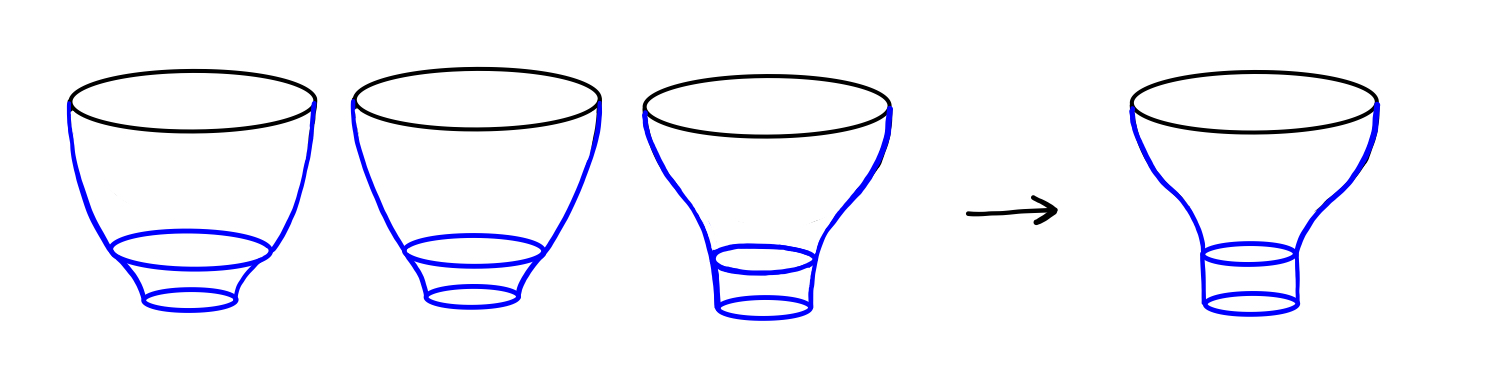} 
\caption{A sequence satisfying the hypotheses of Conjecture~\ref{Conj-NoMin}
whose smooth limit has interior closed minimal surfaces.}\label{fig-cyl-form}
\end{center}
\end{figure}

In contrast Conjecture~\ref{Conj-NoMin} includes hypotheses that are known not to be preserved in the limit.  Imitating examples of Lee-Sormani in \cite{LeeSormani2}, we can construct sequences of $M_j$ satisfying the hypotheses of 
Conjecture~\ref{Conj-NoMin} which converge smoothly to a limit $M_\infty$ that has many interior minimal surfaces.
In fact the region near the boundary can converge to cylinder as seen in Figure~\ref{fig-cyl-form}.    Although no examples have yet been constructed, it is also not clear that limit spaces necessarily have boundaries that are minimal surfaces, so this is not part of the conjecture.

\section{A Compactness Conjecture involving Mass}

A third compactness conjecture has boundaries that are mean convex and minimal.   It first appeared in
\cite{Sormani-scalar} although we have revised it slightly here:

\begin{conj}{\em \bf Scalar Mass Compactness} \label{Conj-Hawking}
Given a sequence of 
three dimensional oriented manifolds $M_j^3$ satisfying 
\be
\vol(M_j) \le V_0\qquad 
\vol(\partial M_j) = A_0 
\qquad \Diam(M_j) \le D_0.
\ee
\be\label{scalar-hawking-bounds}
\Scal_j \ge 0 \qquad \qquad H_{\partial M_j} > 0 \qquad
\qquad m_H(\partial M_j)\le m_0 \qquad
\ee
and either no closed interior minimal surfaces or $\MinA(M_j) \ge A_1>0$,
then a subsequence converges in the intrinsic flat sense
\be
M_{j_k} \VFto M_\infty \textrm{ and } \mass(M_{j_k}) \to \mass(M_\infty)
\ee
where $M_\infty$ is a geodesic metric space which satisfies (\ref{scalar-hawking-bounds})
in some generalized sense.
We might replace Hawking mass with another quasilocal mass in this conjecture.
\end{conj}

\begin{question}
Does this conjecture require the volume bound as a hypothesis or is the upper bound on volume a consequence
of the other hypotheses?
\end{question}

\begin{question}
Can one find a uniform lower bound on the volume depending on the parameters in the hypotheses?   This is
related to a question of Miao and Gromov discussed in Remark~\ref{rmrk-Miao-vol} and \cite{Gromov-Natural}. 
If the conjecture is true, then there is a nonzero limit and so there is a uniform lower bound by the volume convergence.
 \end{question}

LeFloch-Sormani have proven this conjecture in the rotationally
symmetric setting without the volume bound assuming that there are no closed interior minimal
surfaces in \cite{LeFloch-Sormani}.  In fact they proved $H^1_{loc}$ 
convergence to the $\mathcal{F}$ limit and thus the nonnegative scalar curvature
is defined in the $H^1_{loc}$ sense as are the mean curvature and the Hawking mass.

\begin{rmrk}
It is important to note $H^1_{loc}$ convergence 
does not imply $\mathcal{F}$ convergence.  We believe that counter examples can be found (without
scalar curvature bounds) by
adapting Allen and Sormani's examples 
which fail to have $g_j \ge g_0-1/j$ in \cite{AS-relating}.  This would be a nice project for a masters/doctoral student.
See Allen-Bryden \cite{Allen-Bryden-19} for more about Sobolev convergence and $\mathcal{F}$ convergence.
\end{rmrk}

\begin{rmrk}
One possible approach to this conjecture is to assume that $M_j$ can be extended to an asymptotically flat manifold and
that $\partial M_j$ is a level set of the inverse mean curvature flow so that one can apply Huisken-Ilmanen's
equations relating Hawking mass to the geometric properties of the level sets \cite{Huisken-Ilmanen-01}.   
Indeed a key step in the LeFloch-Sormani proof of this conjecture in the rotationally symmetric setting was
to use the monotonicity of the Hawking mass.  It would be interesting to see if
we can prove this compactness theorem in the setting where $M_j$ are covered by a smooth
inverse mean curvature flow from an inner minimal boundary to an outer boundary.   Such regions might
be studied using techniques developed by Brian Allen in 
\cite{Allen-Sobolev}.  
\end{rmrk}

\begin{rmrk}
There is nothing particularly special about the Hawking mass
in Conjecture~\ref{Conj-Hawking}.   We could assume Brown-York Mass is uniformly bounded above
instead or just that the sequence of manifolds have extensions with a uniform upper bound on the ADM mass
of the extension.   The LeFloch-Sormani proof works for all these masses as well because they are
all uniformly bounded at the same time in the rotationally symmetric setting.   
\end{rmrk}

\begin{rmrk}
Examples of Lee-Sormani in \cite{LeeSormani2} demonstrate that, 
under the hypotheses of Conjecture~\ref{Conj-Hawking}, even if a sequence converges smoothly, the
limit space may have a cylindrical region near the inner boundary. so the limit space
may contain closed interior minimal surfaces.   Examples of Mantoulidis and Schoen 
\cite{Mantoulidis-Schoen-neck} 
and of
Cabrera Pacheco, Cederbaum, McCormick, and Miao \cite{CCMM-neck} should also be kept in mind.   
See also the
important work of Corvino in \cite{Corvino-deform}.
\end{rmrk}

\begin{rmrk}
It is possible that Gromov's symmetrization and band approach described in \cite{Gromov-metric} might be useful
towards proving this conjecture and others within this paper.  
\end{rmrk}

Further discussion of this conjecture appears later.


\section{Geometric Stability of the Scalar Torus Rigidity Theorem}   \label{Sect-Torus}

\begin{figure}[h]
\begin{center}
\includegraphics[width=0.8\textwidth]{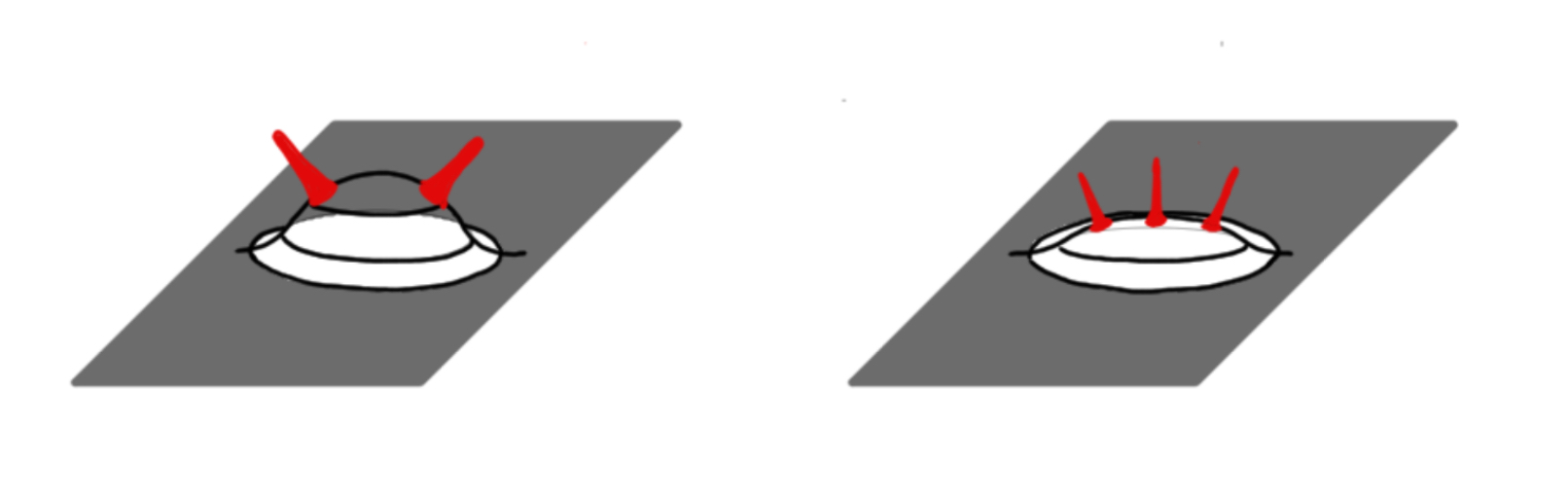} 
\caption{A sequence of Riemannian manifolds satisfying the hypotheses of
Conjecture~\ref{Torus-Stab} which has no GH limit.}\label{fig-tori-wells}
\end{center}
\end{figure}

Let us begin with the geometric stability of the scalar torus rigidity theorem where our manifolds are still homeomorphic to three dimensional tori but we only have {\em almost nonnegative scalar curvature}: $\Scal_j \ge -1/j$.   
See Figure~\ref{fig-tori-wells}.

\begin{conj} {Geometric Stability of Scalar Torus Rigidity [Gromov-Sormani]} \label{Torus-Stab}
Suppose $M_j^3$ are three dimensional Riemannian manifolds 
homeomorphic to tori with 
\be
\Scal_j \ge -1/j \textrm{ and } \MinA(M_j^3) \ge \alpha >0
\ee
satisfying
\be
\Diam(M_j) \le D \textrm{ and } \Vol(M_j) \le V 
\ee
then a subsequence $M_{j_k}^3 \VFto M_\infty$ where $M_\infty$ is isometric to a flat torus.
\end{conj}

If we replace the scalar curvature condition with Ricci curvature, this conjecture was
proven by Colding (see Lemma 3.11 in \cite{Colding-volume}).   He did not need the $\MinA$ condition, just a uniform lower bound on
volume to prevent collapsing.   Note that a sequence of collapsing flat tori, $M_j={\mathbb S}^1 \times {\mathbb S}^1 \times {\mathbb S}^1$ with
\be
g_j=(1/j)^2 dx^2 + dy^2 +dz^2 
\ee
has  
\be
\MinA(M_j)=(2\pi/j)(2\pi)\to 0 \textrm{ and } \Vol(M_j)=(2\pi/j)(2\pi)(2\pi)\to 0,
\ee
so they $\mathcal{VF}$ converge to the $0$ space.

Gromov vaguely stated the conjecture for scalar curvature in \cite{Gromov-Dirac}.  There he proved that if a sequence of $M_j^3$
satisfying the hypotheses converged in the $C^0$ sense to the limit $M_\infty$, then $M_\infty$ is isometric to a flat torus.
Bamler proved the same result in \cite{Bamler-Gromov} using Ricci flow.    Without $C^0$ convergence, the upper bound on diameter and volume  must be added to avoid convergence
to a cylinder or Euclidean space.  With the upper bound on diameter and volume, Wenger's Compactness Theorem implies there is a $\mathcal{F}$ limit \cite{Wenger-compactness}.   Due to the possibility of the existence of increasingly many
increasingly thin wells, as in Figure~\ref{fig-tori-wells}, we cannot hope for GH convergence.

\begin{figure}[h]
\begin{center}
\includegraphics[width=0.8\textwidth]{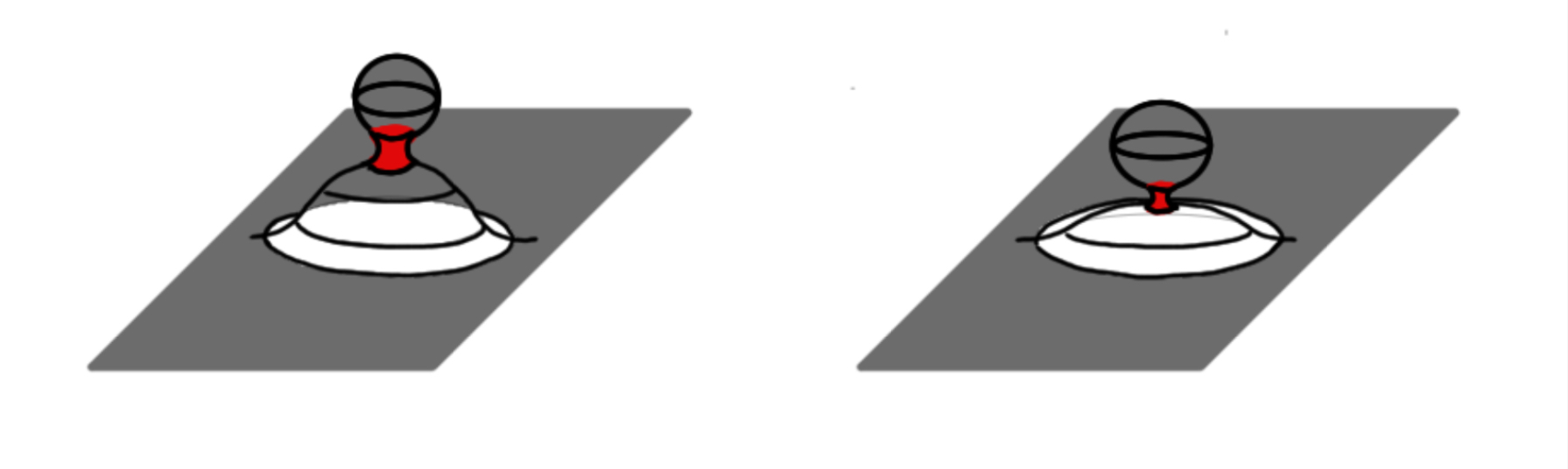} 
\caption{A sequence of Riemannian manifolds satisfying all the hypotheses of
Conjecture~\ref{Torus-Stab} except the $\MinA$ hypothesis which forms bubbles
and fails to converge to a torus.}\label{fig-tori-bubbles}
\end{center}
\end{figure}

The  $\MinA$ condition was added to the hypothesis by Sormani when stating the conjecture in \cite{Sormani-scalar}
in light of examples with bubbling which will appear in work of Basilio-Sormani \cite{BS-tori}.  See Figure~\ref{fig-tori-bubbles}.     We briefly describe the construction here.  Take a standard flat torus with a fixed ball of radius $r_0<\pi/4$.   On this ball we deform the metric tensor to 
\be
g_j=h_j(r)^2 dr^2+f_j(r)^2 g_{\mathbb S}^2 \textrm{ for }r \le \pi/8
\ee
so that $h_j(r)=1$ and $f_j(r)=\sin(\sqrt{K_j}r)/\sqrt{K_j}$ for $r \le \pi/16$.   If we allow $K_j \to 0$ we
can obtain such $g_j$ which converge smoothly to a flat torus.   Then $\Scal_j \to 0$ as well.   Once we have
a region of constant sectional curvature on the manifold we can add a well (cf. \cite{LeeSormani1}) or
increasingly many wells attached at a point as in Figure~\ref{fig-tori-wells}.   We can also attach a tunnel as in \cite{BDS-sewing}, and at the
far side of the tunnel attach it to a sphere of constant curvature $K>K_j$.  In the examples Basilio-Sormani construct in \cite{BS-tori}, the tunnels becomes thinner and thinner as $K_j \to 0$, which causes bubbling: the $\mathcal{VF}$ and
$mGH$ limit is proven to be a torus attached to a sphere of constant curvature $K$.   The $\MinA$ condition prevents the formation of bubbles in this way because minimal surfaces form inside the tunnels with increasingly small area.   

\begin{rmrk}
One also needs the $\MinA$ condition to avoid examples converging to tori which have pulled regions identified to a point which will appear in work of Basilio-Sormani \cite{BS-tori}.   Those examples involving sewing techniques just like those in \cite{BDS-sewing} and \cite{BS-seq} applied to the constant curvature sections described above.  
\end{rmrk}

\begin{rmrk}
Kazaras and Sormani believe that it may also be possible to construct a sequence satisfying all the hypotheses other than the $\MinA_F$ hypothesis that converges to a torus with the taxi distance.   This example would contradict the conjectured conclusion that the limit space has Euclidean tangent cones.
\end{rmrk}

\begin{quest}
Do the hypotheses of Conjecture~\ref{Torus-Stab} imply the existence of a uniform lower bound on volume?  This would follow from the conjecture, because a counter example would imply there was a collapsing sequence
converging to the $0$ space.  It would be interesting to investigate this question directly.
\end{quest}

The conjecture has been proven in the special case where the manifold is a warped product by
Allen, Hernandez-Vazquez, Parise, Payne, and Wang in \cite{AHPPW}.   That is they assume the metric tensor is
of the form
\be
a_j^2(z)\, dx^2 + b_j^2(z) \,dy^2 + dz^2 \quad \textrm{ or } \quad dx^2+dy^2+f_j^2(x,y)\, dz^2
\ee
on ${\mathbb{S}}^1\times {\mathbb{S}}^1\times {\mathbb{S}}^1$.  In fact they are able to prove uniform convergence, mGH convergence and $\mathcal{VF}$ convergence to a flat torus because neither wells not bubbles can form under these strong hypotheses.   

\begin{rmrk}
It would be interesting to prove Conjecture~\ref{Torus-Stab} for locally warped product metric tensors of the form $g_j$ where $g_j= dx^2 +dy^2+dz^2$ away from a collection of balls and $g_j=h_j(r)^2dr^2+f_j(r)^2 g_{{\mathbb S}^2}$ within each ball.    Such sequences can develop increasingly many increasingly thin wells (one in each ball) using the construction described above applied to increasingly many increasingly small balls.   Techniques that might help prove the conjecture in this setting appear in the work of Lee-Sormani \cite{LeeSormani1} and LeFloch-Sormani \cite{LeFloch-Sormani}.   Another approach would be to prove VADB convergence and then apply \cite{VADB} of Allen-Perales-Sormani.
\end{rmrk}

Cabrera Pacheco, Ketterer, and Perales have proven Conjecture~\ref{Torus-Stab} in \cite{CKP-torus-graph} assuming the manifold is a graph over a flat torus, $(M, g_0)$.  That is they assume
\be
g_j = g_0 + df_j^2 \textrm{ where } f_j: \, {\mathbb{T}}^n \to \mathbb{R}
\ee    
with additional hypotheses on $f_j$ related to the hypotheses in earlier work of Lam \cite{Lam-thesis}
and Huang-Lee \cite{HuangLee15}.
It is worth noting that they first announced this paper at IAS with stronger hypotheses on $f_j$ than they required in the published version.  They weakened the hypotheses after an in depth discussion with Gromov and Huang.   In the
final published version, they prove VADB and Sobolev convergence for a subsequence when the base torus, $({{\mathbb{T}}^n},g_{0})$, is fixed.  They then prove $\mathcal{VF}$ convergence of a subsequence in the setting where the base tori are allowed to change using a diagonalization argument.  

\begin{rmrk}
It would be interesting to prove Conjecture~\ref{Torus-Stab} assuming that $M_j$ are tori with metric tensors 
\be
g_j = g_0 + \sum_{k=1}^N dx_j^{k\,2} \textrm{ where } x^k_j: \, {\mathbb{T}}^n \to \mathbb{R}.
\ee   
We would need only prove $\Vol_j(M)\to \Vol_0(M)$ and then we could apply VADB convergence to obtain the $\mathcal{VF}$ convergence.   We could then prove $\mathcal{VF}$ convergence of a subsequence in the setting where the base tori are allowed to change using the same diagonalization argument
that appears in \cite{CKP-torus-graph}.   Note that this setting is quite flexible allowing for the existence of increasingly many wells.    It also seems to allow for the existence of bubbles if we abandon the $\MinA$ condition.   
\end{rmrk}

\begin{rmrk}   
It would be interesting to prove Conjecture~\ref{Torus-Stab} assuming that $M_j$ are conformally flat tori with metric tensors 
\be
g_j = f_j^2\, g_0  \textrm{ where } f_j: \, {\mathbb{T}}^n \to \mathbb{R}.
\ee   
In \cite{Allen-conf-torus}, Allen has proven the conjecture replacing the $\MinA$ hypothesis with a hypothesis that
prevents the concentration of volumes.    Thus we need only prove such a concentration of volume causes the existence of minimal surfaces of increasingly small area.
\end{rmrk}

\begin{rmrk}
It would be interesting to explore if Daniel Stern's new approach in \cite{Stern-Scalar-Harmonic} to proving the torus rigidity theorem  can be applied to prove Conjecture~\ref{Torus-Stab} perhaps with some additional hypotheses.
Recall that he proves his theorem by constructing harmonic maps to $\mathbb{S}^1$ and studying their level sets.
If we assume these harmonic maps are uniformly Lipschitz for a sequence of manifolds, $M_j \VFto M_\infty$, then we might apply
the Sormani's Arzela-Ascoli Theorem in \cite{S-ArzAsc} to study the limits of these harmonic maps defined on $M_\infty$.   Alternatively we might consider harmonic functions on $M_\infty$ and study the level sets as limits of level sets in $M_j$ adapting methods of Portegies-Sormani \cite{PS-properties} and Jauregui-Lee \cite{Jauregui-Lee-isoper}.  One of the great advantages of $\mathcal{VF}$ convergence is control of level sets of Lipschitz functions.
\end{rmrk}

\subsection{Consequences of Torus Stability}

\begin{rmrk}
Recall that Portegies proved $M_j \VFto M_\infty$ implies the spectrum semi-converges \cite{Portegies-evalue}.   
So Conjecture~\ref{Torus-Stab} implies that that the spectrum of $M_j$ semi-converge to the spectrum of a flat torus.
It is impossible here to survey all the results concerning the spectra of tori and flat tori and the first few eigenvalues on tori.
It would be worth studying the consequences of Portegies' result and possibly trying to prove these consequences 
directly.
\end{rmrk}

\begin{rmrk}
We believe that the filling volume of a sphere in Euclidean space with the restricted distance from Euclidean space, is the volume of the Euclidean ball:
$$
\Vol(B(0,R)) = \FillVol(\partial B(0,R)):= \inf\left\{ \mass(N) \, | \, (\partial N, d_N)=(\partial B(0,R), d_{B(0,R)})\right\},
$$
where the infimum is taken over all
integral current spaces $N$ whose boundary with the restricted distance from $N$ is isometric to the boundary of the Euclidean ball with the restricted distance from Euclidean space.   On any flat torus there is a radius $r>0$ such that
$B(p,r)$ with the restricted distance from the torus is isometric to $B(0,r)$ with the restricted distance from Euclidean space.  Consider $M_j \VFto M_\infty$ where $M_\infty$ is a torus.  As mentioned above if $p_j \in M_j$ do not disappear then the filling volumes and volumes converge.  Thus we expect
\be
|\FillVol(\partial B(p_j,r))-\Vol(B(p_j, r))| \to 0.
\ee
This would even be true for disappearing $p_j$ since both would converge to $0$.  Is there any possible way of proving this directly without proving Conjecture~\ref{Torus-Stab} first?
\end{rmrk}

\subsection{Compactness and Torus Stability}

If we prove the MinA Scalar Compactness Conjecture, then the $M_j$ in Conjecture~\ref{Torus-Stab}
would have a converging subsequence $M_{j_k}  \VFto M_\infty$.   We might then be able to prove that $M_\infty$
is a flat torus after obtaining the convergence.  After all there are versions of the Scalar Torus Rigidity Theorem
proven for spaces satisfying weak notions of nonnegative scalar curvature.  Recall that Gromov has proven rigidity
for spaces satisfying the Prism Property in \cite{Gromov-Dirac}.  For this reason the MinA Scalar Compactness Conjecture was formulated at IAS requiring that the limit space satisfy the Prism Property and in particular the limit space 
should have a notion of dihedral angle and mean convexity.     

\begin{rmrk} \label{rmrk-taxi}
It is important to note that if we do not require some notion of angle in our limit space then there are natural counter examples to torus rigidity.  Recall that the Lott-Villani and Sturm $CD(n,0)$ notion of generalized nonnegative Ricci curvature on metric measure spaces defined using optimal transport in \cite{Lott-Villani}\cite{Sturm} includes the taxicab space with Hausdorff measure:
\be
({\mathbb{R}}^3, d_{taxi}, {\mathcal H}^3) \textrm{ where } d_{taxi}((x_1,x_2,x_3),(y_1,y_2,y_3))=\sum_{i=1}^3 |x_i-y_i|.
\ee
This is the universal cover of tori that satisfy the $CD(n,0)$-notion of nonnegative Ricci curvature that are not isometric to flat tori.   Note that these taxi tori have zero scalar curvature in the sense defined by (\ref{Scal-vol-lim}).  While one can easily imagine sequences of three dimensional Riemannian manifolds with $\Scal \ge 0$ that
converge to a taxi torus in the GH sense, these sequences are created using thinner and thinner tunnels around tighter and tighter 3D lattices, so their $\MinA(M_j)\to 0$.    Thus Conjecture~\ref{Conj-MinA} and Conjecture~\ref{Torus-Stab} do not include such sequences.
\end{rmrk}

\begin{rmrk}
Ambrosio-Gigli-Savare added an infinitesimally Euclidean condition to define $RCD(n,K)$ spaces in \cite{Ambrosio-Gigli-Savare} and proved that the $mGH$ limits of noncollapsing sequences of manifolds with nonnegative Ricci curvature lie in this class.  Gigli-Rigoni have proven that $RCD(0,n)$ spaces that are tori are isometric to flat tori \cite{Gigli-Rigoni-18}.  
Thus one way to prove Conjecture~\ref{Torus-Stab} is to first prove the limit space is an $RCD(0,n)$ space.  We know that not all limits obtained in Conjecture~\ref{Conj-MinA} are $RCD(0,n)$ spaces because any Riemannian manifold with negative Ricci curvature at a single point fails this condition.  Thus we would need to prove it is $RCD(0,n)$ using the fact that the limit is a torus not just a limit.   
\end{rmrk}

\begin{rmrk}
Lee and LeFloch have defined and studied nonnegative scalar curvature defined as a distribution on a Riemannian manifold $(M,g)$ with metric tensor $g\in C^0 \cap W^{1,p}(M)$ in \cite{Lee-LeFloch}.  Jiang-Sheng-Zhang have proven torus rigidity for in this setting when $p\in [n,\infty]$  in \cite{JSZ21}.   They first prove that $(M,g)$ is an 
$RCD(0,n)$ space and then apply Gigli-Rigoni's Theorem that $RCD(0,n)$ spaces that are tori are isometric to flat tori \cite{Gigli-Rigoni-18}.     It would be interesting to explore under what additional hypotheses can we guarantee that a
sequence of manifolds converges to a Riemannian manifold with a
$C^0 \cap W^{1,p}(M)$ metric that has distributional $\Scal \ge 0$.   Would VADB convergence suffice?    See related work by Allen-Bryden in \cite{Allen-Bryden-19}.
\end{rmrk}

\begin{rmrk}
Burkhard-Guim has proven torus rigidity for spaces satisfying her weak notion of scalar curvature defined using Ricci flow \cite{B-G-GAFA}.   Under what hypotheses on the sequence would we be guaranteed to have a limit space with a $C^0$ metric in her class?   She has
investigated $C^0$ convergence but it would also be worthwhile to investigate when her property persists under VADB convergence (which allows for the formation of wells).
\end{rmrk}

\begin{rmrk}
Lee-Naber-Neumayer have proven a compactness theorem for sequences of manifolds with nonnegative scalar curvature
satisfying entropy bounds in \cite{Lee-Naber-Neumayer-dp}.   It is unclear if their entropy bounds are satisfied by sequences developing increasingly thin wells.  Their proof involves Ricci flow and it would be interesting to relate their work to that of Burkhardt-Guim and Jiang-Sheng-Zhang mentioned above.   They apply this theorem to obtain a torus stability theorem in all dimensions.   They explain that their notion of $d_p$-convergence is not geometric: more precisely they obtain Gromov-Hausdorff convergence of the
Riemannian manifolds viewed as metric spaces using the distance
\be
d_p(x,y)=\sup \{ \, |f(x)-f(y)| \, | \, \int_M |\nabla f|^p \, dvol_g \le 1,\,\, f\in W^{1,p}\cap C^0_{loc}\}
\ee
rather than the usual distance $d_g$ defined by the infimum of lengths of curves.   They  explain is not geometric.   In dimension 3, it would be interesting to explore how their notion of convergence relates to $VADB$ and $\mathcal{VF}$ convergence.    They do discuss the relationship to various examples of Allen-Sormani in \cite{AS-contrasting} and
\cite{AS-relating} but it is unknown what the $d_p$ limits are of the examples of Basilio, Dodziuk, Kazaras, Sinaei, and Sormani appearing in \cite{BDS-sewing} \cite{BKS-no-geod} \cite{BS-seq} \cite{BS-tori} \cite{Sinaei-Sormani}.
\end{rmrk}

\begin{rmrk}
Keep in mind that Conjecture~\ref{Torus-Stab} does not immediately follow if we prove the MinA Scalar Compactness Conjecture where the limits are in a class of spaces that satisfies the Scalar Torus Rigidity Theorem.   We must also prove that the limit is homeomorphic to a torus.  This becomes an interesting question on its own.  If $M_j$ satisfy the 
hypotheses of Conjecture~\ref{Torus-Stab} and $M_j \VFto M_\infty$, is $M_\infty$ homeomorphic to a torus?   An investigation into the change in topology under $\mathcal{F}$ convergence with no assumption on curvature
appears in work of Sinaei-Sormani \cite{Sinaei-Sormani}.
\end{rmrk}

\section{Geometric Stability of Scalar Prism Rigidity [Gromov-Li]} \label{Sect-Prism}

The geometric stability of the Scalar Prism Rigidity Theorem was formulated at the IAS Emerging Topics Workshop
in 2018 in discussions with Misha Gromov and Chao Li.  As with torus almost rigidity conjecture above, the scalar curvature is assumed to be almost negative, $\Scal_j \ge -1/j$.   We add diameter and volume bounds from above to keep the sequence bounded and a bound on the area of the boundary allows us to use Wenger's Compactness Theorem to find a $\mathcal{F}$ converging subsequence.  

To avoid collapse and other counter examples, we require a lower bound on the area of minimal surfaces with free boundary.   Recall that $\Sigma \subset M$ is
a free boundary minimal surface if it is either a closed minimal surface or it has a boundary, $\partial \Sigma \subset \partial M$ and local variations sliding along the boundary are minimizing.   So in particular $\Sigma$ intersects the boundary at right angles.    We define
\be \label{MinA_F}
\MinA_F(M^3) =\min\{ \Area(\Sigma)\, | \, \Sigma \textrm{ is a free boundary minimal surface in }M^3 \}.
\ee
Note that we also include $\Sigma$ with no boundary in this minimum.

\begin{conj}{Geometric Stability of Scalar Prism Rigidity [Gromov-Li]} \label{Prism-Stab}
If $P_j \subset M_j^3$ are mean convex domains with 6 faces, 8 corners, and 12 edges 
such that the dihedral angles, $\theta_p$, at any $p$ lying on 
an edge of $P^3$ have $\theta_p \le \pi/2$  and 
\be \label{prism-sm}
\Scal_j \ge -1/j \textrm{ and } \MinA_F(M_j^3) \ge \alpha >0
\ee
satisfying
\be\label{prism-Wenger}
\Diam(M_j) \le D \textrm{ and } \Vol(M_j) \le V \textrm{ and } \Area(\partial M_j) \le A
\ee
then a subsequence $P_{j_k}^3 \VFto P_\infty$
where $P_\infty$ is isometric to a rectangular prism in Euclidean space.
\end{conj}

Note that by Wenger's Compactness Theorem and (\ref{prism-Wenger}) we know a subsequence $\mathcal{F}$ converges possibly to the zero space.  First one needs to show the limit is not the zero space using the fact that it
is 3 dimensional and (\ref{prism-sm}).   In dimension 4, one can probably adapt the example of 
Lee-Naber-Neumayer \cite{Lee-Naber-Neumayer-dp} to obtain a counter example satisfying all the other hypotheses.  
Without $\MinA_F$ one can find collapsing flat prisms $[0,1/j] \times [0,1]\times [0,1]$.  Kazaras and Sormani believe that it may also be able to construct a sequence satisfying all the hypotheses other than the $\MinA_F$ with volume bounded below which converges to a cube with the taxi distance or something similar.

At this time Conjecture~\ref{Prism-Stab} is completely open with no partial solutions.   
However the techniques applied to complete the special cases of the Geometric Stability of the Scalar Torus Conjecture should apply equally well to prove the corresponding special cases of this conjecture.    Graduate students might be asked to prove the conjecture in the warped product setting using the techniques developed by Allen-Hernandez-Parise-Payne-Wang in \cite{AHPPW}
or in the graph setting using techniques developed by Cabrera Pacheco, Ketterer, and Perales have proven Conjecture~\ref{Torus-Stab} in \cite{CKP-torus-graph}.  A graduate student might also consider the case where the sequence is known to converge in the VADB sense with boundary as in the work of Allen-Perales \cite{Allen-Perales-VADB} and verify that the limit is flat.   

In light of Gromov's work in \cite{Gromov-Dirac} we have considered whether it is possible that Conjecture~\ref{Torus-Stab} and Conjecture~\ref{Prism-Stab} are equivalent, and whether one might then apply Torus Stability to prove Prism Stability and finally use both to help prove the compactness in Conjecture~\ref{Conj-MinA}.     We discuss this further in the next three remarks which arise from conversations with Misha Gromov, Pengzi Miao, and Chao Li.   These three remarks are part of an overarching vision and should not be viewed as solvable by a doctoral student.   

\begin{figure}[h]
\begin{center}
\includegraphics[width=0.8\textwidth]{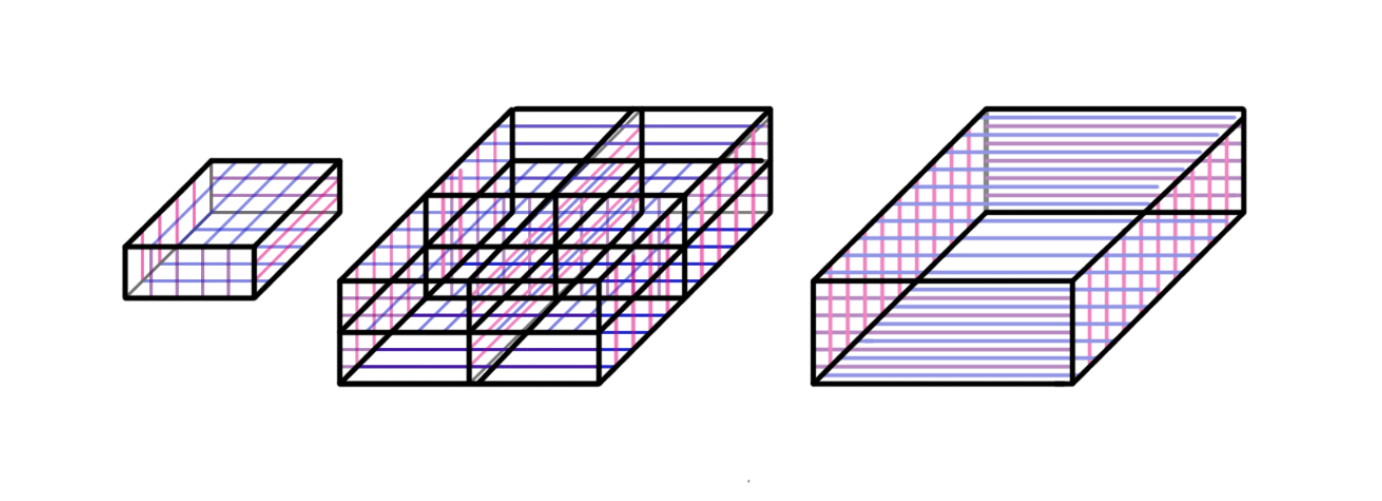} 
\caption{Starting with a prism on the left, Gromov reflects forward, up, and to the right.  Then he attaches the
opposite sides and smoothens the metric to obtain a torus on the right as explained in
Remark~\ref{rmrk-prism2torus}.}\label{fig-prism2torus}
\end{center}
\end{figure}

\begin{rmrk}  \label{rmrk-prism2torus}
Suppose one has a sequence of
prisms satisfying the hypotheses of Conjecture~\ref{Prism-Stab}.  Miao and Li suggest that we might imitate Gromov's technique reflecting each prism across its faces and glueing it together to form a torus as in Figure~\ref{fig-prism2torus}.   We would have a sequence of tori that can be smoothed retaining the lower bound on scalar curvature.  We could check if the tori satisfy the hypotheses of Conjecture~\ref{Torus-Stab}.  This is one reason Li decided to have the $\MinA_F$ condition in Conjecture~\ref{Prism-Stab}.    If the Torus Stability Conjecture holds then we would have a subsequence which
converges to a flat torus.  By \cite{S-ArzAsc}, sufficiently small balls about points in the tori which do not disappear converge to
Euclidean balls.  This would guarantee that the $\mathcal{F}$ limit of of the prisms did not disappear.  Then we might again apply this to prove that the $\mathcal{F}$ limit of the prisms is flat in the sense that small balls about points in the limit are Euclidean balls.  More difficult would be describing the boundary of the limit.
\end{rmrk}

\begin{figure}[h]
\begin{center}
\includegraphics[width=0.8\textwidth]{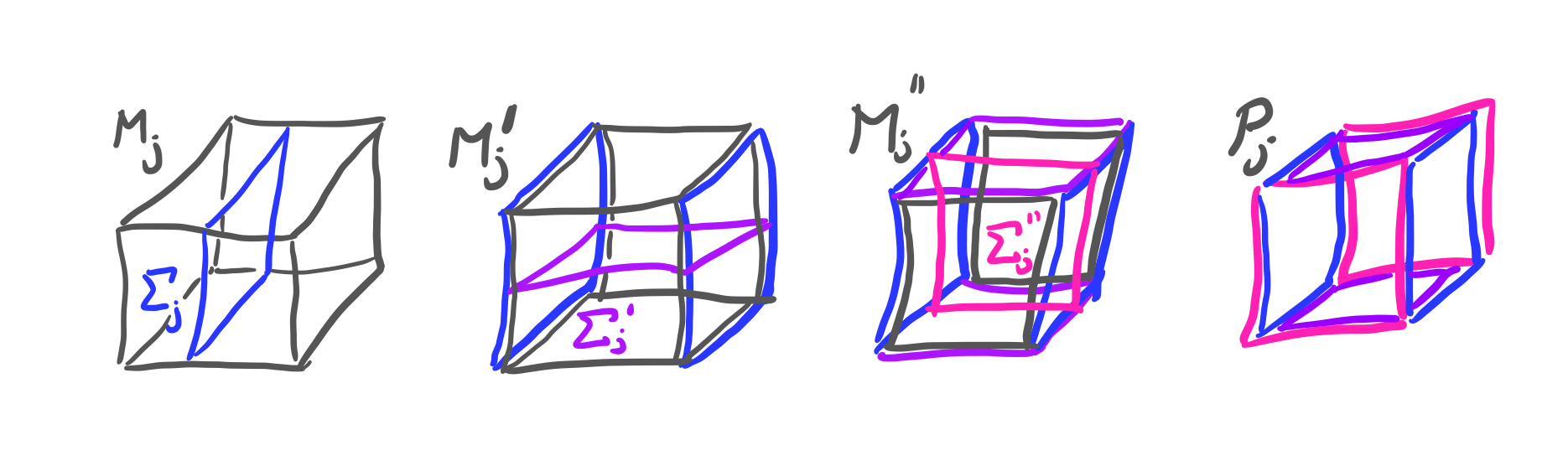} 
\caption{Starting with a torus on the left, we cut along a stable minimal surface, and then along free boundary minimal surfaces to form a prism as described in Remark~\ref{rmrk-torus2prism}.}\label{fig-torus2prism}
\end{center}
\end{figure}

\begin{rmrk} \label{rmrk-torus2prism}
Suppose one has a sequence of tori, $M_j$, satisfying the hypotheses of Conjecture~\ref{Torus-Stab}.   Miao and Li note that we could find a sequence of stable minimal surfaces, $\Sigma_j \subset M_j$, as in Figure~\ref{fig-torus2prism}, and then cut open the $M_j$ along the $\Sigma_j$ to get a manifold, $M_j'$. diffeomorphic to a torus crossed with an interval.   Next one takes a stable free boundary minimal surface $\Sigma'_j\subset M_j'$ and cut the manifold open again to set a manifold $M_j''$ diffeomorphic to a circle crossed with a square region.   Finally we find another stable free boundary minimal surface $\Sigma''_j\subset M_j''$ and cut the manifold open again to get a prism $P_j$ diffeomorphic to a cube.
We could check if the $P_j$ satisfy the hypotheses of Conjecture~\ref{Prism-Stab}, and assuming this conjecture holds, find a subsequence which converges to a flat prism.   By \cite{S-ArzAsc}, sufficiently small balls about points in the $P_j$ which do not disappear converge to
Euclidean balls.  This would guarantee that the $\mathcal{F}$ limit of of the original tori did not disappear.  Then we might again apply this to prove that the $\mathcal{F}$ limit of the tori is flat in the sense that small balls about points in the limit are Euclidean balls.  More difficult would be proving that the limit is a torus.
\end{rmrk}

\begin{rmrk}
Note that Prism Stability would be very useful as a step towards proving that the limits obtained in
Conjecture~\ref{Conj-MinA} satisfy Prism Rigidity.   We would need to first prove the limits are not zero, and that
one can even define a prism on the limit space.   Then we might imitate theorems from \cite{S-ArzAsc} or \cite{Jauregui-Lee-isoper}, to find regions in the $M_j$ converging to the prism in the limit space.   We would then wish to show those regions could be chosen to be prisms satisfying the hypotheses of Conjecture~\ref{Prism-Stab}.  It is quite possible we might need to adjust the conjecture to include more general regions that $\mathcal{VF}$ converge to prisms in order to achieve this.  
\end{rmrk}

\section{Geometric Stability of Scalar Sphere Rigidity [Marques-Neves]}\label{Sect-Sphere}
The following conjecture was proposed by Marques and Neves soon after the IAS Emerging Topics meeting in 2018.   They've graciously allowed us to include the conjecture here:

\begin{conj}{Geometric Stability of Scalar Sphere Rigidity [Marques-Neves]} \label{Sphere-Stab}
Suppose $M_j^3$ homeomorphic to spheres with
\be
\Scal_j \ge 6-1/j \textrm{ and } \MinA(M_j^3) \ge 4\pi-1/j
\ee
satisfying
\be \label{Sphere-DM}
\Diam(M_j) \le D \textrm{ and } \Vol(M_j) \le V 
\ee
then $M_{j}^3 \VFto {\mathbb S}^3$
where ${\mathbb S}^3$ is the standard round three sphere.
\end{conj}

\begin{rmrk}
Note that in this conjecture Marques and Neves force the width to be almost maximal using the fact that
the width is achieved by a minimal surface and thus:
\be
\Width(M_j^3) \ge \MinA(M_j^3) \ge 4\pi.
\ee
In ${\mathbb S}^3$ we have equality above and $\Scal_j=6$.  If $MinA({\mathbb S}^3, g) > 4\pi/3$ and $\Scal_g \geq 6$, then $({\mathbb S}^3,g)$ has no stable minimal spheres. In that case, it follows from \cite{MN-Duke} that $\MinA({\mathbb S}^3,g)=Width({\mathbb S}^3,g)$. They proved in that paper that if $R_g\geq 6$ and $A(g)\geq 4\pi$, then $g$ has constant sectional curvature one.
\end{rmrk}

\begin{rmrk}
Note that the volume and diameter bounds in the conjecture guarantee a subsequence converges in the intrinsic flat sense by Wenger's Compactness Theorem \cite{Wenger-compactness}.   The goal is to prove the limit must be a round sphere and that the volume converges to the volume of the sphere.    It would be an interesting first step to prove that under the hypothesis of this theorem we have
\be
\Vol(M_j) \to \Vol({\mathbb S}^3).
\ee
\end{rmrk}

\begin{figure}[h]
\begin{center}
\includegraphics[width=0.8\textwidth]{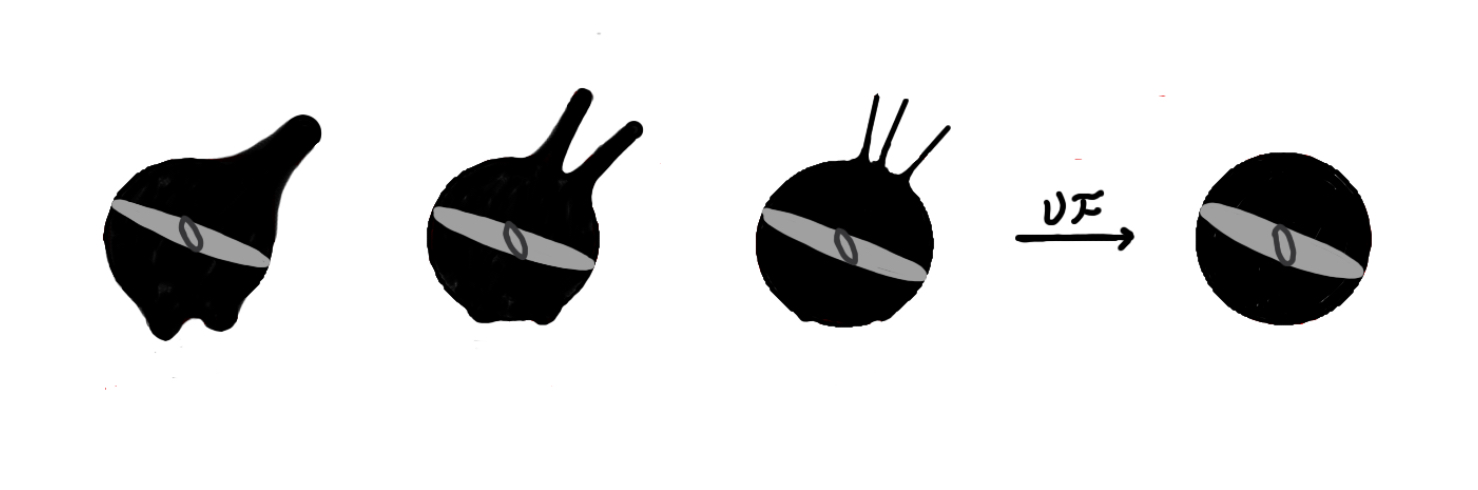} 
\caption{Here is an example of a sequence of $M_j^3$ with unstable minimal surfaces $\Sigma_j$
such that $\Area(\Sigma_j)=\MinA(M_j)=\Width(M_j^3)$ such that $M_j^3 \VFto {\mathbb S}^3$.  However, it
is unknown if $\Scal_j \ge 6-1/j$ at the base of the wells.}\label{fig-maybe2sphere}
\end{center}
\end{figure}

\begin{rmrk}
Note that there are not yet examples precisely constructed related to this conjecture.  We believe that someone should perhaps be able to construct sequences of spheres with increasingly many increasingly thin wells which satisfy the hypothesis of Conjecture~\ref{Sphere-Stab} as in Figure~\ref{fig-maybe2sphere}.  Such a sequence would have no GH limit but would converge in the VADB and $\mathcal{VF}$ sense to a standard sphere.   However it is not completely clear that a well construction has been precisely defined with $\Scal \ge 6-1/j$ rather than just $\Scal >0$.   Note that the long thin part of the well and the tip easily have large scalar curvature, but the part where the well attaches to the sphere would need to be calculated precisely.
\end{rmrk}
 
\begin{rmrk}
Bamler and Maximo have proven this conjecture with the additional hypothesis that $M_j$ have positive sectional curvature in \cite{Bamler-Maximo}.  In fact with such strong additional hypothesis they are able to prove $C^0$ convergence.   Their proof builds on work of Brendle and Schoen \cite{Brendle-Schoen} and can also be applied to prove almost rigidity or geometric stability of the $\mathbb{RP}^3$ Rigidity Theorem of Bray-Brendle-Eichmair-Neves \cite{BBEN10}.   The sectional curvature condition is quite strong, and without it we expect there are counter examples to $C^0$ convergence as  in Figure~\ref{fig-maybe2sphere}.
\end{rmrk}
   
\begin{rmrk}
Daniel Stern has proven the scalar sphere rigidity theorem using his level sets of harmonic maps approach in \cite{Stern-Scalar-Harmonic}.  Perhaps his approach might be applied to prove Conjecture~\ref{Sphere-Stab}.  We might test this
with some additional hypotheses.   As mentioned earlier, his approach might be combined with work of
Jauregui, Lee, Portegies, and Sormani in \cite{Jauregui-Lee-isoper}\cite{S-ArzAsc} \cite{PS-properties}.
\end{rmrk}

\begin{rmrk}
As a test for this conjecture, we might consider $M_j$ as in the conjecture and assume
$M_j \VADBto M_\infty$ and try to prove that $M_\infty$ is a standard sphere.   
\end{rmrk}

\begin{rmrk}
Recall that Portegies proved $M_j \VFto M_\infty$ implies the spectrum semi-converges \cite{Portegies-evalue}.   
So Conjecture~\ref{Sphere-Stab} implies that that the spectrum of $M_j$ semi-converge to the spectrum of a sphere.
It is impossible here to survey all the results concerning the spectra of spheres and their eigenfunctions.
It would be worth studying the consequences of Portegies' result and possibly trying to prove these consequences 
directly.  It is even possible that one might use the one of the spectral rigidity theorems for a sphere as a step
towards proving Conjecture~\ref{Sphere-Stab} in its full generality or at least a special case.
\end{rmrk}

\begin{rmrk}
It is believed that the filling volume of a standard round three sphere with the restricted distance from Euclidean space, is the volume of a four dimensional hemisphere:
\be
\Vol({\mathbb{S}}^4_+) = \FillVol({\mathbb{S}}^3):= \inf\{ \mass(N) \, | \, (\partial N, d_N)
=({\mathbb{S}}^3, d_{{\mathbb{S}}^3})\},
\ee
where here $\mass(N)$ is the Ambrosio-Kirchheim mass and the infimum is taken over all
integral current spaces $N$ whose boundary with the restricted distance for $N$ is isometric to the standard sphere with its standard distance.   Observe that the four dimensional Euclidean ball is not a filling for the sphere because its boundary with the restricted distance is not isometric to the sphere.  In \cite{PS-properties}, Portegies-Sormani proved the filling volumes converge, so if one has $M_j$ satisfying the hypotheses of Conjecture~\ref{Sphere-Stab}, and
the conjecture is true, we expect
\be
\FillVol(M_j) \to \Vol({\mathbb{S}}^4_+).
\ee
It would be interesting to check if one is able to prove this directly without proving the conjecture first.
\end{rmrk}

\begin{rmrk}
Gromov has stated a different Spherical Stability Conjecture in \cite{Gromov-Natural} based on Llarull's Scalar Rigidity Theorem \cite{Llarull}.   Llarull proved that if there is degree one distance nonincreasing map from $M^n$ to ${\mathbb{S}}^n$
and $\Scal(M^n) \ge n(n-1)$ then $M^n$ is isometric to ${\mathbb{S}}^n$.  Gromov has suggested that if there
are degree one distance nonincreasing maps from $M_j^n$ to ${\mathbb{S}}^n$
and $\Scal(M_j^n) \ge n(n-1)-(1/j)$ then $M^n \Fto {\mathbb{S}}^n$.   We expect a $\MinA$ condition might be
required to avoid bubbling when $n=3$ and possibly something stronger in higher dimensions.
\end{rmrk}

\section{Geometric Stability of Zero Mass Rigidity [Lee-Sormani]}\label{Sect-Mass}

The final conjecture concerns the geometric stability of the various zero mass rigidity theorems.   This conjecture was first proposed by Dan A. Lee and Christina Sormani in \cite{LeeSormani1} where it was proven in the spherically symmetric setting.     

\begin{conj} {Geometric Stability of Zero Mass Rigidity [Lee-Sormani]}\label{Mass-Stab}
Suppose $M_j^3$ are asymptotically flat with no closed interior minimal surfaces and possibly minimal boundary
with $\Scal \ge 0$ everywhere and $m_{ADM}(M^3)\le 1/j$.   Suppose $\Omega_j\subset M_j^3$ are
regions inside $M_j$ containing $\partial M_j$ whose outer boundaries 
$\Sigma_j=\partial \Omega_j \setminus \Omega_j$ have constant mean curvature 
and
\be \label{DVA-mass}
\Diam(\Omega_j) \le D \textrm{ and }  \Area(\Sigma_j)=A_0
\ee
then $\Omega_j \VFto \Omega_0$ where $\Omega_0$ is a ball in Euclidean space such that
$\Area(\partial \Omega_0)=A_0$.
\end{conj}

See the final section of \cite{LeeSormani1} for a variety of restatements of this conjecture with different requirements on the regions $\Omega$.  The regions might have their boundaries within uniformly asymptotically flat regions of $M_j$ or
they might have positive Gauss curvature.  Note also that we might consider compact regions $\Omega_j$ satisfying these hypotheses whose Hawking masses or Brown-York masses converge to $0$ without assuming they have asymptotically flat extensions $M_j$ whose ADM masses converge to $0$.  

\begin{rmrk}\label{mass-vol-conv}
Note that proving that $\Vol(\Omega_j) \to \Vol(\Omega_0)$ under the hypotheses of Conjecture~\ref{Mass-Stab} would be an important result on its own.   Keep in mind that without the diameter bounds preventing increasingly deep wells the volumes can diverge to infinity.  Without the hypothesis that there are no closed interior minimal surfaces, bubbles could form preventing the volume from converging as well.   Both these phenomena can occur even in the spherically symmetric setting.   It is possible that Huisken's Isoperimetric mass and work of Fan-Shi-Tam 
in \cite{FST09} might help, but one needs more control in the interior.  See Remark~\ref{rmrk-Miao-vol}.
\end{rmrk}

In order to prove a geometric stability or almost rigidity theorem like Conjecture~\ref{Mass-Stab}, we must look for estimates on the mass of the manifold or the quasilocal mass of the region that depend upon the geometry of the interior. 
The Schoen-Yau proof of the Zero Mass Rigidity Theorem is a proof by contradiction, so it cannot be applied to control the $M_j$.    In the special cases of this conjecture that have been proven,  we do have such formulas for the mass.  In Lee and Sormani's proof in the spherically symmetric setting they use the monotonicity of the Hawking mass of the spheres which increases to the ADM mass.   With the Hawking mass of all spheres close to zero, they were able to control the warping function strongly away from a possible central well whose volume converges to 0.

\begin{figure}[h]
\begin{center}
\includegraphics[width=0.6\textwidth]{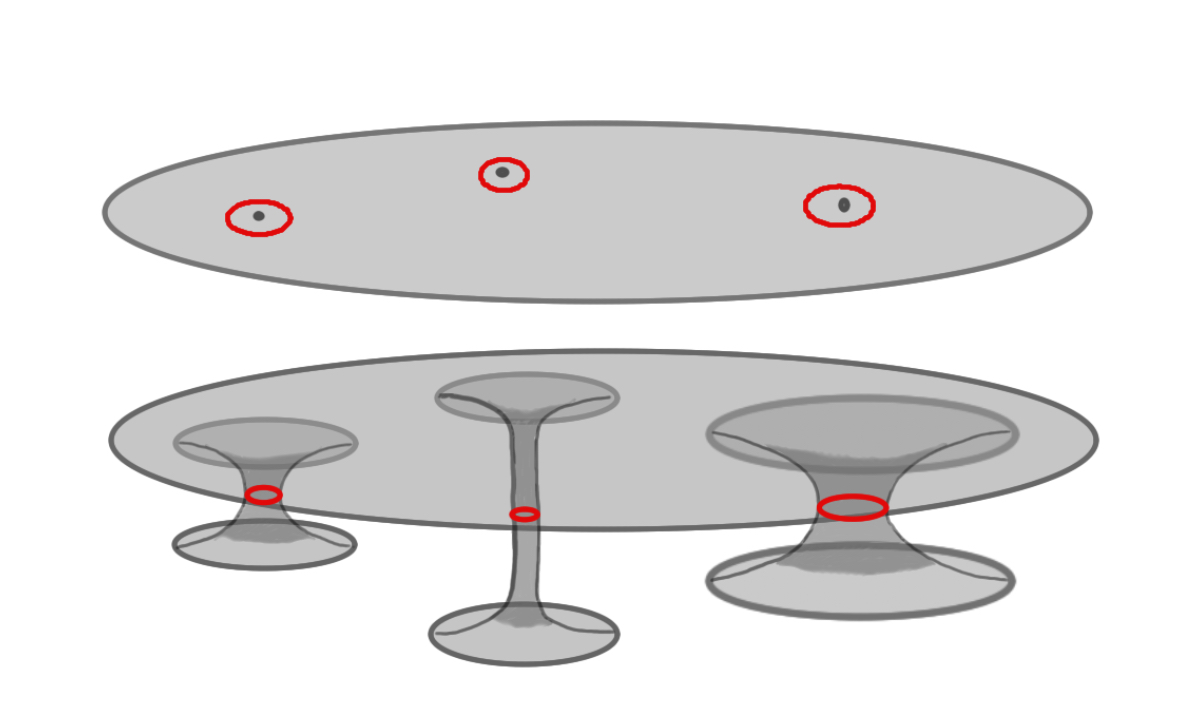} 
\caption{This geometrostatic manifold $M'$ has three poles $p_i$ with necks whose depth and area  depend on $a_i$ and $b_i$.  The region $M \subset M'$ lying above the minimal surfaces
 satisfies the hypothesis of Conjecture~\ref{Mass-Stab} and is $\mathcal{VF}$ close to ${\mathbb E}^3.$}\label{fig-geometrostatic}
\end{center}
\end{figure}

\begin{rmrk}
In \cite{Sormani-Stavrov}, Sormani and Stavrov have proven Conjecture~\ref{Mass-Stab}
for regions in geometrostatic manifolds $M'_j$ of the form
\be
({\mathbb{R}}^3\setminus \{p^j_1,....,p^j_N\}, g_j) 
\ee
where $g_j$ is conformally flat
\be
g_j=\left( 1 + \sum_{i=1}^{N_j} \frac{a^j_i}{\rho_i}\right)^2\left( 1 + \sum_{i=1}^N \frac{b^j_i}{\rho_i}\right)^2 g_0
\ee
with $a_i>0$, $b_i>0$, and $\rho_i(x)=|x-p_i|$.  In such manifolds 
\be
m_{ADM}(M_j)=\sum_{i=1}^{N_j} (a^j_i+b^j_i).
\ee
So $m_{ADM}(M_j)\to 0$ immediately implies that regions which are a definite distance away from the poles converge smoothly
to the Euclidean metric $g_0$.  However, each pole is a new asymptotically flat end.  See Figure~\ref{fig-geometrostatic}.
To find $M_j \subset M_j'$
which satisfy the hypotheses of Conjecture~\ref{Mass-Stab}, Sormani and Stavrov cut along minimal surfaces
whose locations were previously unknown.  They locate the surfaces in annuli by assuming the $a_i$ and $b_i$ are small relative to $\min|p_i-p_j|$.   They then find domains in the
$M_j$ which are biLipschitz close to Euclidean domains.   They complete the proof by estimating the $\mathcal{F}$
distance using the work of Lakzian-Sormani \cite{Lakzian-Sormani}.   It would be interesting to study where the minimal
surfaces are when the $a_i$ and $b_i$ are large relative to $\min|p_i-p_j|$.   
\end{rmrk}

\begin{rmrk}
We believe the techniques of Sormani-Stavrov in \cite{Sormani-Stavrov}, of Benjamin-Stavrov in \cite{BA18},
of Benjamin-McDermott-Stavrov
in \cite{BMS}, and more recent work of Stavrov announced at VWRS
could be applied to prove Conjecture~\ref{Conj-Hawking} for regions in geometrostatic  manifolds with uniform upper bounds on their ADM masses and uniform bounds on the locations of the poles $p^j_i$.   We would first choose a subsequence such that 
\be
p^j_i\to p^\infty_i \textrm{ and } a^j_i\to a^\infty_i \textrm{ and } b^j_i\to b^\infty_i 
\ee
and then we would have to devise an hypothesis that guarantees we can locate the minimal surfaces for a subsequence and complete the proof.   We would recommend this as doctoral dissertation project and suggest that students interested in this project contact Sormani and Stavrov to ensure they have reserved it for themselves.   
\end{rmrk}

\begin{rmrk}
Conjecture~\ref{Mass-Stab} has also been proven in the graph setting under a variety of hypotheses by 
Huang-Lee-Sormani \cite{HLS}, Allen-Perales \cite{Allen-Perales-VADB}, and Huang-Lee-Perales \cite{HLP}. The
original approach in \cite{HLS} in which Huang-Lee-Sormani first proved $\Omega_j \VFto \Omega_\infty$ and then proved $\Omega_\infty$ is a Euclidean ball was very intriguing but the endplay was not completed rigorously.   It would be very interesting to complete their proof as originally outlined but this would require a deeper analysis of Ambrosio-Kirchheim theory.  Those who are interested should communicate with Perales  as she has studied this problem closely.  The proof of the conjecture in the graph setting was eventually completed correctly by Allen, Huang, Lee, and Perales in \cite{Allen-Perales-VADB} and \cite{HLP} by first proving VADB convergence using Allen, Perales, and Sormani's work in \cite{VADB}.   Note that the techniques in these papers might be applied to prove the Hawking Compactness Conjecture~\ref{Conj-Hawking} for geometrostatic  manifolds.   This could be quite challenging as the VADB convergence might not work and so we do not recommend it as a doctoral dissertation project.
\end{rmrk}

\begin{rmrk}
In light of Zero Local Mass Rigidity Theorem of Shi-Tam in \cite{Shi-Tam-JDG02}, we expect Conjecture~\ref{Mass-Stab} should hold for regions $\Omega_j$ whose Brown-York mass converges to $0$.  Initial work in this direction has been completed in the graph setting by Alaee, Cabrera Pacheco, and McCormick in \cite{AlCaMc-quasi}.  
\end{rmrk}

\begin{rmrk}
There are also asymptotically hyperbolic versions of this conjecture with $\Scal \ge -6$ proven by Sakovich and Sormani in the spherically symmetric setting \cite{Sakovich-Sormani-hyp} and announced work by Cabrera Pacheco and Perales in the graph setting \cite{CabreraPeralesHyp}.   It would be intriguing to see if there might be a asymptotically hyperbolic manifolds with $\Scal \ge -6$ similar to the geometrostatic manifolds studied by Sormani-Stavrov in \cite{Sormani-Stavrov}.   If so, we would suggest trying to prove the asymptotically hyperbolic version of this conjecture in that setting as well.  See also Chao Li's paper
\cite{Li-dihedral-hyp}.
\end{rmrk}

\begin{rmrk}  
There are partial results towards proving Conjecture~\ref{Mass-Stab} in special settings which do not quite achieve $\mathcal{VF}$ convergence but do obtain Sobolev bounds on the metric tensors.   
Allen controls regions $\Omega_j$ that 
are covered by smooth inverse mean curvature flow with bounds on their Hawking masses in \cite{Allen-Sobolev}.   
Bryden controls regions in axisymmetric manifolds in \cite{Bryden20}.   For a survey of earlier results see
\cite{LeeSormani1}. In these papers the metric tensors are shown to converge with various levels of regularity and some also prove volume convergence.   We believe these results might now be improved to achieve VADB and thus $\mathcal{VF}$ convergence using the more recent methods developed by Allen-Perales-Sormani in \cite{VADB}.
However it is not yet clear how we can obtain lower bounds on distances in these settings.
\end{rmrk}

\begin{rmrk}
As with the other stability conjectures, we may consider $M_j$ satisfying the hypotheses and assume the $\Omega_j$ converge in a stronger sense to a limit $M_\infty$ and then prove $M_\infty$ is Euclidean space.
For example, if we assume $C^0$ convergence to a smooth $M_\infty$ we can use Gromov's Prism Property as in \cite{Gromov-Dirac} to conclude that $M_\infty$ has $\Scal\ge 0$.  Jauregui and Lee have proven in \cite{Jauregui-Lee-isoper} that under $C^0$ convergence 
\be \label{JL-mass}
m_{ADM}(M_\infty) \le \liminf_{j \to \infty} m_{ADM}(M_j)
\ee
and then we can apply Zero Mass Rigidity for smooth manifolds to complete the proof.  If one assumes less regularity on the limit space, but assumes the convergence is $C^0$ and $H^1_{loc}$, then we could apply low regularity zero mass rigidity theorems of Lee and LeFloch  \cite{Lee-LeFloch} and Jiang-Sheng-Zhang  in \cite{JSZ21} to prove
$M_\infty$ is Euclidean space.  
\end{rmrk}

\subsection{Using Compactness Theorems}

Note that Wenger's Compactness Theorem and the hypotheses of Conjecture~\ref{Mass-Stab} guarantee a subsequence of the $\Omega_j$ $\mathcal{F}$ converge to some $\Omega_\infty$ which might be the $0$ space. We would "only" need to show that
$\Omega_\infty$ is always a Euclidean ball and thus the original sequence converges as required in the conjecture.   If we assume the $M_j$ are uniformly asymptotically flat in a strong enough way for $R\ge R_0$ that their asymptotically flat ends converge to an asymptotically flat end $E_\infty$ and we require that $\partial \Omega_j$ lie in this asymptotically flat regions, we might hope to prove
$\partial \Omega_j$ cannot converge to $0$ space.  This approach was suggested in Sormani's survey article \cite{Sormani-scalar} and attempted in the graph setting by Huang-Lee-Sormani in \cite{HLS}.

We must ensure that extrinsic distances between points in $\partial \Omega_j$ do not all converge to $0$ which requires knowledge of geodesics which pass through $\Omega_j$.  Basilio-Sormani have constructed examples in  \cite{BS-seq} where entire regions in manifolds with positive scalar curvature disappear in the limit because all distances between points in those regions converge to $0$.  These examples contain increasingly many tiny tunnels and so they have many closed interior minimal surfaces.   To prove $\partial \Omega_j \Fto \partial \Omega_\infty\neq 0$ we could try controlling distances from below using the fact that there are no closed interior minimal surfaces in the hypotheses of Conjecture~\ref{Mass-Stab}.  This is related to Remark~\ref{rmrk-geods}.  Keep in mind also that in light of the examples of Lee-Naber-Neumayer we need the hypothesis that we are in dimension three as well
\cite{Lee-Naber-Neumayer-dp}.

Alternatively we could require as a hypothesis that $\partial \Omega_j$ lie deep within 
the uniformly asymptotically flat ends where $R\ge 2R_0$ so that 
when their distances are strongly estimated by
\be 
d_{M_j}(p,q) \ge d_{E_j|\partial E_j}(p,q) \quad \forall p,q\in \partial \Omega_j
\ee
where
\be
d_{E_j|\partial E_j}(p,q)=\min\{ d_{E_j}(p,q), \min\{d_{E_j}(p,x)+d_{E_j}(q,y)\, |\, x,y \in \partial E_j\} .
\ee
If so,  we believe it might be possible to prove that
\be
\FillVol(\partial \Omega_j, d_{M_j}) \ge \FillVol(\partial \Omega_j, d_{E_j|\partial E_j}) \to 
\FillVol(\partial \Omega_\infty, d_{E_\infty|\partial E_\infty}) >0
\ee
because 
\be
\partial \Omega_\infty=R^{-1}(2R_0) \textrm{ and } E_\infty={\mathbb E}^3 \setminus B(0, R_0).
\ee
This would require new theorems about filling volumes
but intuitively it seems true.  If we are able to prove this then we can try gluing $\Omega_\infty$ to $E_\infty$ to obtain a limit space $M_\infty$ that we hope to prove is Euclidean using a weak version of the Zero Mass Rigidity Theorem.   The Arzela-Ascoli Theorems in \cite{S-ArzAsc} can be useful for gluing.   Assuming all of this works we have $\mathcal{F}$
convergence to a region $\Omega_\infty\subset M_\infty$.

If someone has proven volumes converge as discussed in Remark~\ref{mass-vol-conv}, then we have $\mathcal{VF}$ convergence.  Jauregui and Lee have
proven semicontinuity of (\ref{JL-mass}) for $\mathcal{VF}$ limits in \cite{Jauregui-Lee-VF} using isoperimetric methods.   Thus we would have $m_{ADM}(M_\infty)=0$.   We could then complete this vaguely outlined proof of Conjecture~\ref{Mass-Stab} by proving $M_\infty$ has some low regularity notion of nonnegative scalar curvature and then proving the Zero Mass Rigidity Theorem for such limit spaces.  In particular proving the Prism Property holds on the limit space and proving the Prism Stability Conjecture might complete the proof.

Alternatively we could use one of our stronger conjectured compactness theorems to prove $\Omega_j\VFto \Omega_\infty$ and that $\Omega_\infty\neq 0$ and satisfies the Prism Property.   We would still need to use many of these outlined steps even in this setting.  Note also that the concerns mentioned in this outline are also concerns when proving the compactness conjectures and in those settings we would not be able to use the trick  where we assume 
$\partial \Omega_j$ lie deep within the uniformly asymptotically flat ends unless we extend the manifolds to create such ends.

It is essential that we be very cautious not to simply solve one part of this outline or a few parts and then claim one of the conjectures is true following the approach suggested here.   Anyone proving these theorems needs to complete the proof carefully themselves checking every step or piecing together the proof carefully citing different papers and verifying the hypotheses match.  This outline has been presented to suggest interesting problems to study and to warn people away from jumping to conclusions too quickly.  Significant further study must be completed to better understand the properties of filling volumes as defined in \cite{PS-properties} and the properties of integral current spaces.

\subsection{Possible Approaches to a Quantitative Proof}

To directly prove Conjecture~\ref{Mass-Stab} we must estimate $d_\mathcal{F}(\Omega_j, \Omega_0)$
and $|\Vol(\Omega_j)-\Vol(\Omega_0)|$.   This approach has been successful in the special cases outlined above.  Here we suggest a few possible approaches towards proving this conjecture in full generality in a direct quantitative way.  If anyone attempts one of these proofs and needs to add an hypothesis to do so then they will have proven an interesting result as long as the hypothesis is natural.   Or perhaps someone might find a counter example, in which case the hypotheses of our conjectures would have to change.

\begin{rmrk}
Huisken has suggested that we apply Huisken and Ilmanen's monotonicity of the Hawking mass under their weak inverse mean curvature flow \cite{Huisken-Ilmanen-01}.  Their flow covers some regions smoothly and those regions might be controlled using the techniques developed by Allen in \cite{Allen-Sobolev}.   Their flow skips over some regions which Huisken has suggested might be shown to have total volume converging to $0$ using ordinary mean curvature flow
back inwards.  The volume estimates of Fan-Shi-Tam  in \cite{FST09} might be useful.
We might then use Lakzian-Sormani's estimates on the intrinsic flat distance as in \cite{Lakzian-Sormani}
to complete the proof.   We expect this would be quite challenging to implement but would be well worth investigating even if we only achieve partial results.  It would even be interesting just to prove the volume convergence perhaps in this way.   
\end{rmrk}

\begin{rmrk}
We must of course suggest that it is possible that an investigation into Witten's proof of the Zero Mass Rigidity
for Spin manifolds in \cite{Witten-PMT} would lead to a partial result.  Work in this direction has been completed
by Lee-LeFloch \cite{Lee-LeFloch}.  There was some concern in the discussions at IAS that the Spin approach might not lead to $\mathcal{VF}$ convergence but a different type of convergence.  The biggest concern was that we could not devise a way to relate Spin to the existence of minimal surfaces.    
\end{rmrk}

\begin{rmrk}
Yu Li has proven Zero Mass Rigidity using Ricci Flow in \cite{Li-RicciflowPMT}.   It is possible his proof can be investigated
in the setting where the ADM mass is small.   Note that in three dimensions he applies Ricci flow with surgery.   It might be helpful to apply the work of Lakzian in \cite{Lakzian-Ricci}
which demonstrates that manifolds which flow through a neck pinch singularity 
in the sense of Angenant-Caputo-Knopf \cite{ACK-neck} are evolving continuously with respect to $\mathcal{VF}$ convergence.
\end{rmrk}

\begin{rmrk} 
In \cite{BKKS-PMT} Bray-Kazaras-Khuri-Stern have proven zero mass rigidity 
using harmonic maps to circles.   Some of these authors are investigating what can 
be determined about the
harmonic maps when the mass is small in this setting.
\end{rmrk}

\section{Additional Related Conjectures and Thoughts}

In light of our discussions throughout it seems essential to investigate the following questions:

\begin{question}
What happens to the Width under $\mathcal{VF}$ and VADB convergence?   Width varies continuously under Lipschitz convergence
as can be seen for example by examining the proofs in Colding and Minicozzi's paper \cite{CoMin-width} which  demonstrates how Width evolves under Ricci flow.   
\end{question}  

\begin{question}
What happens to Hamilton's Ricci flow under $\mathcal{VF}$ or VADB convergence of the initial manifolds?   Simon proved $C^0$ convergence of $C^2$ metric tensors $g_j \to g_\infty$ with $\Scal_j \ge \kappa$ implies smoother convergence of the flows $g_{j,t}\to g_{\infty,t}$ in \cite{Simon02}.   This would be particularly interesting in light
of the work of Bamler in \cite{Bamler-Gromov} and Burkhardt-Guim in \cite{B-G-GAFA}.
\end{question}

\begin{question}
What happens to heat kernels under $\mathcal{VF}$ or VADB  convergence?  It is known how they converge under Lipschitz convergence and weaker settings as proven in work of Yu Ding \cite{Ding-heat}.
\end{question}   

\begin{question}
What happens to harmonic maps under $\mathcal{VF}$ or VADB convergence?   This would be particularly interesting to study in light of work of Daniel Stern \cite{Stern-Scalar-Harmonic}.
\end{question}   

We might also consider additional geometric stability conjectures: 

\begin{question}
Can we formulate and prove a $\mathcal{VF}$ geometric stability conjecture for
Hemispherical Scalar Rigidity Theorem of Eichmair \cite{Eichmair-Hemi09}?
\end{question}

\begin{question}
Can we formulate and prove a $\mathcal{VF}$ geometric stability conjecture for
the Cover Splitting Scalar Rigidity Theorem of Bray, Brendle, and Neves \cite{BBN10}?
\end{question}

\begin{question}
Can we formulate and prove a $\mathcal{VF}$ geometric stability conjecture for
the ${\mathbb{RP}}^3$
Scalar Rigidity Theorem by Bray, Brendle, Eichmair and Neves \cite{BBEN10}?
\end{question}

\begin{question}
Can we formulate and prove a $\mathcal{VF}$ geometric stability conjecture for
Hyperbolic Scalar Rigidity Theorem of Nu\~nez \cite{N-JGA-13}?
\end{question}

For each one, it would be natural to weaken the hypothesis on scalar curvature, preserve the other hypotheses, require uniform upper bounds on volume and diameter, and a uniform lower bound on $\MinA$.
However we must be cautious with each such conjecture, particularly with regard to the preservation of the
hypotheses under $\mathcal{VF}$ convergence.  See for example \cite{Sormani-scalar}.   It would be worthwhile to test hypotheses of these other possible conjectures by proving them in special cases like those discussed throughout this paper. 

In addition to the conjectures and questions posted in this article, there are more concerning about $\mathcal{F}$ and $\mathcal{VF}$ convergence without scalar curvature bounds posed by Gromov, Perales, Portegies, Nu\~nez Zimbr\'on, Sinaei, Sormani, and Wenger that appear in \cite{Gromov-Dirac}\cite{Gromov-Plateau}\cite{S-ArzAsc}\cite{PS-properties}\cite{PN}, \cite{Sinaei-Sormani}\cite{Sormani-scalar}.   Progress on these questions can be checked on my website:
$$
\textrm{https://sites.google.com/site/intrinsicflatconvergence/}
$$
which includes a regularly updated bibliography of all papers concerning intrinsic flat convergence
including their abstracts with comments.

\newpage
\section{A final review}

Let us quickly summarize everything:
\be
M_j \LIPto M_\infty \quad \implies \quad M_j \VFto M_\infty 
\ee
was proven by Sormani-Wenger in \cite{SW-JDG} and more constructively by
Lakzian-Sormani in \cite{Lakzian-Sormani}.  In \cite{VADB}. Allen-Perales-Sormani
proved:
\be
M_j \VADBto M_\infty \quad \implies \quad M_j \VFto M_\infty.
\ee
Finally there is the obvious fact that
\be
M_j \VFto M_\infty \quad \implies \quad M_j \Fto M_\infty
\ee
which further implies that 
\begin{itemize}
\item volumes/masses converge \cite{SW-JDG},
\item boundaries converge \cite{SW-JDG},
\item diameters subconverge \cite{S-ArzAsc},
\item balls converge \cite{S-ArzAsc},
\item filling volumes converge \cite{PS-properties},
\item sliced filling volumes converge \cite{PS-properties},
\item unif. bounded Lip. functions to a compact metric space converge \cite{S-ArzAsc},
\item uniformly local isometries converge \cite{S-ArzAsc},
\item there exist
subsets $U_j \subset M_j$ such that $U_j \GHto M_\infty$ \cite{S-ArzAsc}.
\end{itemize}
If one further has $M_j \VFto M_\infty$ then
\begin{itemize}
\item $M_j \mto M_\infty$ \cite{Portegies-evalue}
\item $\limsup_{j\to \infty} \lambda_k(M_j) \le \lambda_k(M_\infty)$ \cite{Portegies-evalue}
\item volumes of balls converge \cite{Jauregui-Lee-VF}
\item areas of spheres converge \cite{Jauregui-Lee-VF}
\item various isoperimetric types of regions converge \cite{Jauregui-Lee-VF}
\end{itemize}
If we assume $M_j \VFto M_\infty$ or $M_j \VADBto M_\infty$, we ask the following questions:
\begin{itemize}
\item What happens to $\Width(M_j)$?
\item What happens to the Ricci flow of $M_j$?
\item What happens to the heat kernals of the $M_j$?
\item What happens to harmonic maps $h_j: M_j \to X$?
\end{itemize}
Hopefully these properties will not only provide new insight into methods towards proving our compactness and geometric stability conjectures but also lead to a deeper understanding of three dimensional oriented compact Riemannian manifolds with lower bounds on their scalar curvature.
\newpage

\bibliographystyle{plain}
\bibliography{2021}

\end{document}